\documentclass[reqno]{amsart}

\pdfoutput=1

\usepackage{scrextend}
\deffootnotemark{\textsuperscript{(\thefootnotemark)}}

\usepackage{hyperref}
\usepackage{microtype}

\usepackage[french, english]{babel}
\usepackage[alphabetic]{amsrefs}
\usepackage[pdftex]{color}
\usepackage{mathtools}
\usepackage{amsxtra}
\usepackage{upref}
\usepackage{epstopdf}
\usepackage{graphicx,color}
\usepackage{hyperref}
\usepackage{longtable}
\usepackage{graphicx}
\usepackage{amsmath, amsfonts, amssymb,amscd, amstext, amsthm, 
mathtools, bezier, mathrsfs,enumerate}
\usepackage{url}
\usepackage{float}
\usepackage{caption}
\usepackage{tikz}
\usetikzlibrary{intersections,positioning,patterns,arrows,decorations.markings,matrix}
\usepackage[fit]{truncate}
\usepackage{anyfontsize}
\usepackage{t1enc}
\usepackage{synttree}
\usepackage{pst-plot}
\usepackage{todonotes}
\psset{algebraic}
\usepackage{arydshln,leftidx}

\newcommand{\h}{h_{\text{top}}}

\newcommand{\vr}{\varrho}
\newcommand{\wh}{\widehat}
\newcommand{\vt}{\vartheta}
\newcommand{\un}{\underline}
\newcommand{\ve}{\varepsilon}

\newcommand{\R}{\mathbb{R}}
\newcommand{\C}{\mathbb{C}}
\newcommand{\Q}{\mathbb{Q}}

\newcommand{\N}{\mathbb{N}}
\newcommand{\F}{\mathcal{F}}
\newcommand{\K}{\mathcal{K}}
\newcommand{\cc}{\mathscr{C}}
\newcommand{\rr}{\restriction}
\newcommand{\f}{\widetilde{f}}
\newcommand{\iS}{\mathring{S}}

\theoremstyle{plain}
\newtheorem{thm}{Theorem}
\newtheorem{lem}[thm]{Lemma}
\newtheorem{prop}[thm]{Proposition}
\newtheorem{cor}[thm]{Corollary}
\newtheorem{stm}[thm]{Statement}
\newtheorem{example}{Example}

\theoremstyle{definition}
\newtheorem{definition}[thm]{Definition}

\theoremstyle{remark}
\newtheorem{rem}[thm]{\bf Remark}

\numberwithin{equation}{section}
\numberwithin{thm}{section}

\makeatletter
\newsavebox{\@brx}
\newcommand{\llangle}[1][]{\savebox{\@brx}{\(\m@th{#1\langle}\)}%
  \mathopen{\copy\@brx\kern-0.5\wd\@brx\usebox{\@brx}}}
\newcommand{\rrangle}[1][]{\savebox{\@brx}{\(\m@th{#1\rangle}\)}%
  \mathclose{\copy\@brx\kern-0.5\wd\@brx\usebox{\@brx}}}
\makeatother

\DeclareFontFamily{OML}{rsfs}{\skewchar\font'177}
\DeclareFontShape{OML}{rsfs}{m}{n}{ <5> <6> rsfs5 <7> <8> <9> rsfs7
  <10> <10.95> <12> <14.4> <17.28> <20.74> <24.88> rsfs10 }{}
\DeclareMathAlphabet{\mathfs}{OML}{rsfs}{m}{n}

\newcommand{\intav}[1]{\mathchoice {\mathop{\vrule width 6pt height 3 pt depth  -2.5pt
\kern -8pt \intop}\nolimits_{\kern -6pt#1}} {\mathop{\vrule width
5pt height 3  pt depth -2.6pt \kern -6pt \intop}\nolimits_{#1}}
{\mathop{\vrule width 5pt height 3 pt depth -2.6pt \kern -6pt
\intop}\nolimits_{#1}} {\mathop{\vrule width 5pt height 3 pt depth
-2.6pt \kern -6pt \intop}\nolimits_{#1}}}

\newcommand{\intavl}[1]{\mathchoice {\mathop{\vrule width 6pt height 3 pt depth  -2.5pt
\kern -8pt \intop}\limits_{\kern -6pt#1}} {\mathop{\vrule width 5pt
height 3  pt depth -2.6pt \kern -6pt \intop}\nolimits_{#1}}
{\mathop{\vrule width 5pt height 3 pt depth -2.6pt \kern -6pt
\intop}\nolimits_{#1}} {\mathop{\vrule width 5pt height 3 pt depth
-2.6pt \kern -6pt \intop}\nolimits_{#1}}}

\begin{document}


\title[Kneading Theory for Iteration of Monotonous Functions]{Kneading Theory for Iteration of Monotonous Functions on the Real Line}
\author{Ermerson Araujo}
\author{Alex Zamudio Espinosa}
\address{Departamento de Matem\'atica Aplicada/IME, Universidade Federal Fluminense (UFF), 
Campus do Gragoat\'a, CEP 24210-201, Niter\'oi -- RJ, Brazil}
\email{alexmze@id.uff.br}
\address{Departamento de Matem\'atica/CCET, Universidade Federal do Maranh\~ao (UFMA), 
Cidade Universit\'aria Dom Delgado, CEP 65080-805, S\~ao Lu\'is -- MA, Brazil}
\email{ermerson.araujo@ufma.br}

\date{\today}
\keywords{Kneading theory, topological conjugacy, monotonous maps}
\subjclass[2020]{37B10, 37B40, 37C15, 37E05, 37F05}

\begin{abstract}
We construct a version of kneading theory for families of monotonous functions on the real line. The generality of the setup covers two classical results from Milnor-Thurston's kneading theory: the first one is to dynamically characterise an $l$-modal map by its kneading sequence, the second one is to define the concept of kneading determinant, relate it to topological entropy and use this to construct a certain type of special ``linearazing measure''.
\end{abstract}

\maketitle

\section{Introduction}

One important problem in dynamical systems is whether we can characterize the dynamics of a system if we only know the combinatorics of some special orbits. The classical kneading theory, first developed by Milnor and Thurston, studies this problem in the context of multimodal maps. They managed to get information about the properties of a multimodal map from the knowledge of the combinatorics of the orbits of some special points, called turning points. 

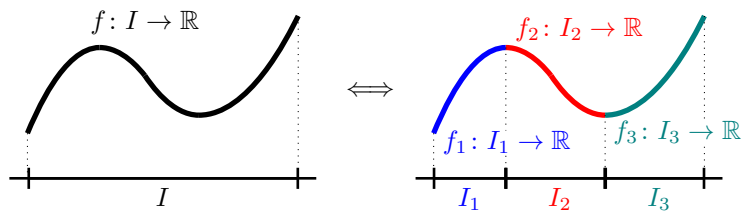
\begin{figure}[H]\centering
\begin{tikzpicture}[scale=1.2]
      \draw[line width=1pt] (0.3,0) -- (3.7,0);
			\draw[|-|][line width=1pt] (0.5,0) -- (3.51,0);
			\draw[line width=2pt] (1.3,1.45) parabola (0.5,0.5);
			\draw[line width=2pt] (1.3,1.45) parabola (1.8,1.14);
			\draw[line width=2pt] (2.4,0.7) parabola (1.8,1.14);
			\draw[line width=2pt] (2.4,0.7) parabola (3.5,1.8);
			\draw[dotted] (0.5,0.5) -- (0.5,0);
			\draw[dotted] (3.5,1.8) -- (3.5,0);
			\fill[black] (1.8,1.5) circle (0.mm) node[above]{$f\colon I\to\R$};
			\fill[black] (2,0) circle (0.mm) node[below]{$I$};
			\fill[black,very thick] (4.3,0.8) circle (0.mm) node[above] {$\Longleftrightarrow$};
			\draw[line width=1pt] (4.8,0) -- (8.2,0);
			\draw[|-|][line width=1pt] (4.99,0) -- (5.81,0);
			\draw[|-|][line width=1pt] (5.79,0) -- (6.91,0);
			\draw[|-|][line width=1pt] (6.892,0) -- (8.01,0);
			\draw[blue,line width=2pt] (5.8,1.45) parabola (5,0.5);
			\draw[red,line width=2pt] (5.8,1.45) parabola (6.3,1.14);
			\draw[red,line width=2pt] (6.9,0.7) parabola (6.3,1.14);
			\draw[teal,line width=2pt] (6.9,0.7) parabola (8,1.8);
			\draw[dotted] (5,0.5) -- (5,0);
			\draw[dotted] (5.8,1.45) -- (5.8,0.53);
			\draw[dotted] (5.8,0.25) -- (5.8,0);			
			\draw[dotted] (6.9,0.7) -- (6.9,0);
			\draw[dotted] (8,1.8) -- (8,0.66);
			\draw[dotted] (8,0.35) -- (8,0);
			\fill[black] (5,0.4) circle (0.mm) node[right]{\textcolor{blue}{$f_1\colon I_1\to\R$}};
			\fill[black] (5.8,1.65) circle (0.mm) node[right]{\textcolor{red}{$f_2\colon I_2\to\R$}};
			\fill[black] (7.7,0.28) circle (0.mm) node[above]{\textcolor{teal}{$f_3\colon I_3\to\R$}};
			\fill[black] (5.4,0) circle (0.mm) node[below]{\textcolor{blue}{$I_1$}};
			\fill[black] (6.4,0) circle (0.mm) node[below]{\textcolor{red}{$I_2$}};
			\fill[black] (7.5,0) circle (0.mm) node[below]{\textcolor{teal}{$I_3$}};
\end{tikzpicture}
\caption{From an $l$-modal map to a system of monotonous functions.}
\label{fig1}
\end{figure}

The purpose of this paper is to extend some of the concepts and results of kneading theory to the context of iteration of monotonous functions on the real line. The basic motivation is that an $l$-modal function (the basic object studied by kneading theory) can be thought as a family of $l+1$ strictly monotonous functions. The study of the dynamic of the $l$-modal function would then be equivalent to the study of the iteration of the monotonous functions;
see Figure \ref{fig1}.

Our belief is that many of the properties and features of kneading theory rely mainly in the fact that those functions are strictly monotonous, nor in the fact that they glue together continuously or even in the fact that their domains only intersect at border points. Previous indications of such a theory are present in other works: in \cite{AMa1990} Alsedà and Mañosas create a version of kneading theory for a special type of function on the circle, then use this to relate the set of possible itineraries with the interval of rotation numbers. 
Using the setting built on \cite{Baillif1999}, Baillif and de Carvalho \cite{BaCa2001} construct a linear model for any continuous piecewise strictly monotone tree map. In \cite{BHV1} Barnsley, Harding and Vince study critical itineraries for a particular type of function which is piecewise continuous and piecewise monotonous and give a formula for its topological entropy, in this paper we prove the same formula using our tools. In \cite{AM} Alsed\`a and Misiurewicz prove that any piecewise monotonous, piecewise continuous function with positive entropy can be semiconjugated with a piecewise affine function with constant slope. We obtain a similar result assuming the function is also piecewise strictly monotonous.

The basic object of study in this paper will be a family of real valued functions $\F=\{f_1,\ldots ,f_n\}$, where each $f_i$ is strictly monotonous with domain a closed interval. The objective will be to extend two classic results from kneading theory. The first one is to dynamically characterise an $l$-modal map by its kneading sequence (as in Chapter II, Section 3 of \cite{dMvS}), the second one is to define the concept of kneading determinant, relate it to topological entropy and use this to construct a certain type of special measure (as in Theorem 6.3 and Section 7 from \cite{MT}).

We would also like to highlight that our approach somehow brings together different aspects of one dimensional dynamics. For example, Theorem \ref{thm:measureLamb} guarantees the existence of a ``linearazing measure'' for an $l$-modal map with positive entropy, but it also gives an stationary measure for an iterated function system of contractions defined on a closed interval. Another example is that Theorem \ref{main:thm:1} serves to characterise $l$-modal maps by its kneading sequences, but in the context of dynamically defined Cantor sets in the real line it basically corresponds to the fact that the topological dynamics of the Cantor set is determined by the transition matrix associated to the Markov partition.

The statement of our theorem characterizing the dynamics of $\F$ using its kneading sequences is the following:

\begin{thm}\label{main:thm:1}
Let $\F=\{f_1, \ldots, f_n\}$ and $\widetilde{\F}=\{\f_1, \ldots, \f_n\}$ be two systems with ${\rm Dom}(f_i)={\rm Dom}(\f_i)$ and such that the orientations of $f_i$ and $\f_i$ are equal, 
for all $1\leq i\leq n$, and assume that both satisfy the separability hypothesis. Then: 
\begin{enumerate}[{\rm (1)}]
\item $\F$ and $\widetilde{\F}$ are combinatorially equivalent if and only if they have the same kneading sequence. 
\item In addition, also assume that if an interval $(x,y)$ is a connected component of $\R\setminus {\rm cl}(\cc_G)$ then $x,\, y \in \cc_G$, and similarly for $\widetilde{\cc}_G$. Then, $\F$ in ${\rm cl}(\cc_G)$ is topologically conjugate to $\widetilde{\F}$ in ${\rm cl}(\cc_G)$ if and only if they have the same kneading sequence.
\end{enumerate}
\end{thm}

For $\F=\{f_1,\ldots ,f_n\}$ we define its topological entropy in the following way: for each positive integer $m$ define $\ell(m)$ as the number of sequences $(i_1,\ldots,i_m)$, with $i_j\in \{1,\ldots,n\}$, such that 
$$\text{int(Dom} (f_{i_1}\circ\cdots\circ f_{i_m}))\neq \emptyset.$$
The growth number of $\F$ is defined as
$s=\limsup_{m\to \infty} \ell(m)^{1/m}.$
The topological entropy of $\F$ is then defined as $h_{\text{top}}(\F)= \log s.$
This definition is motivated by the formula of Misiurewicz and Szlenk (see \cite{MS1980}) for the topological entropy of an $l$-modal map. In Section \ref{entropy} we show that $h_{\text{top}}(\F)$ coincides with the topological entropy of the skew-product associated to $\F$.

For a system $\F$ we also define the number $s_0$, which in some sense measures the growth of individual orbits. For example, if $\F$ comes from an $l$-modal function, then individual orbits do not grow and $s_0=1$. We will see that positive entropy, in the context of $l$-modal maps, corresponds to $s> s_0$ in our context. When this happens we will prove that the kneading determinant $D(t)$ vanishes at $t=1/s$ and it does not vanishes for $|t|<1/s$. This is the content of Theorem \ref{thm:ss0}, which corresponds to our version of Theorem 6.3 in the classical work of Milnor and Thurston \cite{MT}. 

We will use the results mentioned in the previous paragraph to show the existence of a special measure associated to the system $\F$:

\begin{thm}\label{main:thm:2}
Assume that $s>s_0$. Then there exists a $\sigma$-algebra $\mathscr{B}_{\F}$ in $\R$ and a measure $\Lambda:\mathscr{B}_{\F} \to \R$ such that
\[\Lambda(J)=\tfrac{1}{s} \left[\Lambda(f_{a_1}(J))+\cdots+\Lambda(f_{a_n}(J))\right]\]
for all $J\in \mathscr{B}_{\F}$.
\end{thm}

The paper is organized as follows: in Section \ref{definitions} we give basic definitions and notation. In Section \ref{ks:IFS} we treat the problem of dynamically characterizing the topological conjugacy class of a system using kneading sequences. In Section \ref{MTKT} we define kneading matrix, kneading determinant and relate it to entropy, we also construct a special type of measure associated to the system. In Section \ref{entropy} we define entropy and show that the definition coincides with the entropy of the associated skew product. Finally in Section \ref{appexamples} we give some examples and applications.

\section{Basic definitions}\label{definitions}

Throughout this text, $\F$ will always denote a system $\{f_i\colon I_i\to\R\,:\, i=1,\ldots,n\}$ of continuous strictly monotonous functions, where $I_i=[c_i^\ell,c_i^r]$ is a closed interval in the real line $\R$, for $i=1,\ldots, n$. This will be the basic object with which we will be working, sometimes we will refer to such $\F$ just as a ``system''. The points $c_i^\ell, c_i^r$ will play the role of the turning points in the classical work of Milnor and Thurston, even though they do not look like ``turning points'' in this new context we will still call them {\it turning points} in order to maintain a standard nomenclature of the literature.
In general, we assume that each $I_i$ is non-degenerate, i.e. $I_i$ has more than one point. We will need to consider degenerate intervals for some applications, in this case we will say that the system $\F$ is degenerate. Unless otherwise stated, each system in the text is assumed to be non-degenerate.

We denote by $G$ the free group generated by the $n$ symbols 
$a_1, \ldots, a_n$. The identity element in $G$ is denoted by $e$  and $|g|$ stands for the word-length of the element $g\in G$. More precisely, every $g\in G\setminus\{e\}$ can be written in a unique way as $g=x_1x_2\cdots x_m$, where each $x_i$ is either a generator or the inverse of a generator, and $x_{i+1}\neq x_i ^{-1}$. The word-length of $g$ is defined by $|g|=m$, and $|e|=0$. Given $g,\, h$ two elements of $G$, define the distance between them as $|g^{-1}h|$. We say that a sequence $g_1,\, g_2,\ldots,\, g_m$ of elements of $G$ is a geodesic if the distance between $g_i$ and $g_j$ is $|i-j|$. We will also consider the free semigroup generated by 
the symbols $a_1, \ldots, a_n$ and denote it by $S$. We will denote by $S^{-1}$ the set of inverses of elements in $S$.

For $g\in G$ we define $f_g$ following the rules:
\[
f_{a_i}=f_i,\; f_e=\textrm{Id}_{\R}\; \textrm{and}\; f_{gh}=f_h\circ f_g.
\]
Denote by $\textrm{Dom}(f_g)$ the domain of the function $f_g$, notice that it is possible that $\textrm{Dom}(f_g)=\emptyset$. 
The equation $f_{gh}=f_h\circ f_g$ means that 
$\textrm{Dom}(f_g)\cap f_g^{-1}(\textrm{Dom}(f_h))\subset \textrm{Dom}(f_{gh})$ 
and $f_{gh}(x)=f_h(f_g(x))$ for all $x\in \textrm{Dom}(f_g)\cap f_g^{-1}(\textrm{Dom}(f_h))$. 
Moreover, if we have $|gh|=|g|+|h|$ then $\textrm{Dom}(f_g)\cap f_g^{-1}(\textrm{Dom}(f_h))= \textrm{Dom}(f_{gh})$. 

For a function $f:\text{Dom}(f)\to \R$ and a set $J\subset \R$, not necessarily contained in $\text{Dom}(f)$, when we write $f(J)$ we mean $f(J\cap \text{Dom}(f))$. Given $x\in\R$, write $x_m\uparrow x$ to denote a sequence $(x_m)_{m\in\N}\in\R$ 
such that $x_m\leq x_{m+1}\leq x$ and 
$\lim_{m\to\infty}x_m=x$. Similarly, $x_m\downarrow x$ means that
$\lim_{m\to\infty}x_m=x$ and $x\leq x_{m+1}\leq x_m$. Given a set $X$, we denote its interior and closure by ${\rm int}(X)$ and ${\rm cl}(X)$, respectively.

%


\section{Kneading sequences}\label{ks:IFS}

The aim of this section is to show that, in some sense, the kneading sequences characterize the topological conjugacy class of a system. In fact, we will show that, under certain hypothesis, they determine the combinatorial equivalence class of a system. After that, we will add an extra hypothesis to be able to pass from combinatorial equivalence to topological equivalence. Before proceeding, we give the formal definition of the terms we just mentioned. For that consider a system $\F=\{f_1, \ldots, f_n\}$ like above.


\medskip
\noindent
{\sc Address of a point:}  The {\em address} of a point $x\in \R$ is a 
mathematical object capable of determining the relative position of 
$x$ with respect to $I_i$, for each $i=1, \ldots, n$. For each $x$ and $i$ we have five options: $x\in\textrm{int}(I_i)$, 
or $x\notin I_i$ and it is to the left of the interval, or $x\notin I_i$ and 
it is to the right of the interval, or $x=c_i^\ell$, or $x=c_i^r$. We denote this address
by $D(x)$. For this section we do not need to explicitly define $D(x)$, we just need its properties. At the beginning of Section \ref{MTKT} we give a exact definition. 

\medskip
\noindent
{\sc Admissibility:} Let $x\in \R$ and $g\in G$. We say that 
$g$ is {\it admissible} for $x$ if $x\in \textrm{Dom}(f_g)$. 
We denote the set of admissible elements for $x$ by $G_x$ (analogously $S_x$).

\medskip
\noindent
{\sc Itinerary of a point:} Let $x\in \R$. The {\it itinerary} of $x$ is the 
function $T_x$, defined on $G_x$, taking $g\in G_x$ to $D(f_g(x))$.

\medskip
\noindent
{\sc Kneading sequence:} The {\em kneading sequence} of $\F=\{f_1, \ldots, f_n\}$ 
is the family of functions $\K=\{T_c\rr_{S_c}: c \textrm{ is turning point}\}$.

\medskip
\noindent
{\sc Critical orbits:} We call {\it critical $G$--orbit} of $\F=\{f_1, \ldots, f_n\}$
the set $\{f_g(c): g\in G_c\; \textrm{and}\; c\; \textrm{is turning point}\}$.
We denote this set by $\cc_G$. One defines analogously the {\it critical $S$--orbit}, $\cc_S$, and the {\it critical $S^{-1}$--orbit} $\cc_{S^{-1}}$.

\medskip
\noindent
{\sc Invariant set:} Given $K\subset \R$, we say it is {\it $\F$--invariant} if $f_g(K\cap \text{Dom}(f_g))\subset K$ for all $g\in G$. For example, the critical $G$--orbit $\cc_G$ is invariant.

\medskip
\noindent
{\sc The map $\lambda$:} We define $\lambda:\R\times\R\to\{-1,0,1\}$ as
$$
\lambda(x,y)= \left\{
\begin{array}{rl}
 1   & \textrm{ if }  x<y, \\
 0   & \textrm{ if }  x=y, \\
-1   & \textrm{ if }  x>y.
\end{array}
\right.
$$

\medskip
\noindent
{\sc The map $\sigma$:} Let $\widehat{\sigma}:G \to \{-1,1\}$ be the unique homomorphism such that $\widehat{\sigma}(a_i)=1$ if $f_i$ is increasing and $\widehat{\sigma}(a_i)=-1$ if $f_i$ is decreasing. Given $g\in G$ we define $\sigma(f_g)=\widehat{\sigma}(g)$. 

\medskip

It is easy to see that
\[
\lambda(f_g(x),f_g(y))=\sigma(f_g) \lambda(x,y)
\]
whenever $x,y\in\textrm{Dom}(f_g)$.
In the next lemma we show how to use the itinerary of points to determine the order between them.

\begin{lem}\label{lem:lamb}
Let $x,y \in \R$, we have:
\begin{enumerate}[{\rm (a)}]
\item If $S_x=S_y$ and $T_x\rr_{S_x}=T_y\rr_{S_y}$, then $S_z=S_x$ and $T_z\rr_{S_z}=T_x\rr_{S_x}$ for all $z$ between $x$ and $y$.
\item If $D(x)\neq D(y)$ then $\lambda(x,y)$ can be determined from the values of $D(x)$, $D(y)$.
\item If $S_x\neq S_y$, then one can determine $\lambda(x,y)$ from the values of $T_x\rr_{S_x}$, $T_y\rr_{S_y}$ and the orientation of the functions.
\end{enumerate}
\end{lem}

\begin{proof}
(a) Assume $S_x=S_y$ and take $z$ between $x$ and $y$. Given $g\in S_x=S_y$ we have $x,y \in \textrm{Dom}(f_g)$, since $\textrm{Dom}(f_g)$ is connected we conclude that $z\in \textrm{Dom}(f_g)$. Therefore $S_x\subset S_z$. Now, given $g\in S_z$ there is a sequence $f_{g_0}, \ldots, f_{g_j}$ such that $g_0$ is the identity of $S$, $g_m=g_{m-1}a_{i(m)}$, for some generating element $a_{i(m)}$, $g_j=g$ and $z\in \textrm{Dom}(f_{g_m})$ for any $1\leq m\leq j$. If $g\notin S_x=S_y$ then for some $m$ we have $x,y\in \textrm{Dom}(f_{g_{m-1}})$ and $x,y\notin \textrm{Dom}(f_{g_m})$. 
It follows that the interval limited by $f_{g_{m-1}}(x), f_{g_{m-1}}(y)$ is either to the right or left of $I_{i(m)}$ since $D(f_{g_{m-1}}(x))=D(f_{g_{m-1}}(y))$.
But this would imply $z\notin \textrm{Dom}(f_{g_m})$, we conclude that $g\in S_x=S_y$ and therefore $S_z=S_x=S_y$. Finally, given $g\in S_z=S_x=S_y$ we have that $f_g(z)$ is between $f_g(x)$ and $f_g(y)$, then $D(f_g(z))=D(f_g(x))$ and we conclude $T_z\rr_{S_z}=T_x\rr_{S_x}$.

\medskip
\noindent
(b) Since $D(x)\neq D(y)$ then there is a turning point $c$ between $x$ and $y$. The values of $D(x)$, $D(y)$ determine the values of $\lambda(x,c)$, $\lambda(y,c)$ and therefore determines $\lambda(x,y)$.

\medskip
\noindent
(c) Assume $S_x\neq S_y$, then there is $g\in S_x$ such that $y\notin \textrm{Dom}(f_g)$ or $t\in S_y$ such that 
$x\notin\textrm{Dom}(f_t)$. Suppose the first one.
Like item (a) there exists a sequence $f_{g_0}, \ldots, f_{g_j}$ such that $g_0$ is the identity of $S$, $g_m=g_{m-1}a_{i(m)}$, for some generating element $a_{i(m)}$, $g_j=g$ and $x\in \textrm{Dom}(f_{g_m})$ for any $1\leq m\leq j$.
Then, for some $m$ we have $y\in \textrm{Dom}(f_{g_{m-1}})$ and $y\notin \textrm{Dom}(f_{g_m})$ (notice that $y\in \textrm{Dom}(f_{g_0})=\R$). 
This implies that $D(f_{g_{m-1}}(x))\neq D(f_{g_{m-1}}(y))$, by item (b) above we can 
determine the value of $\lambda(f_{g_{m-1}}(x),f_{g_{m-1}}(y))$ from the 
values of $T_x\rr_{S_x}$, $T_y\rr_{S_y}$, finally we have 
$\lambda(x,y)=\sigma(f_{g_{m-1}}) \lambda(f_{g_{m-1}}(x),f_{g_{m-1}}(y))$.
\end{proof}

In the following we will deal with two systems $\F=\{f_1, \ldots, f_n\}$ 
and $\widetilde{\F}=\{\f_1, \ldots, \f_n\}$, to denote the elements related to the 
second system we will use the symbol $\sim$, for example $\widetilde{T}_x$ 
would be the itinerary of the point $x$ associated to $\widetilde{\F}$. 


\begin{definition}
Let $\{f_1, \ldots, f_n\}$ and $\{\f_1, \ldots, \f_n\}$ 
be two systems with $\sigma(f_i)=\sigma(\f_i)$, for all $1\leq i\leq n$.
We say that they are {\it combinatorially equivalent} if there
exists an order preserving bijection $\varphi:\cc_G\to\widetilde{\cc}_G$ such that $\varphi(c_i^{\ell/r})=\widetilde{c}_i^{\ell/r}$, for $i=1,\ldots,n$, $\varphi(\textrm{Dom}(f_g)\cap\cc_G)=\textrm{Dom}(\f_g)\cap\widetilde{\cc}_G$ and
$\varphi\circ f_g (x)= \f_g \circ\varphi(x)$ for all $g\in G$ and $x\in \textrm{Dom}(f_g)\cap\cc_G$. 
\end{definition}

Observe that the necessary condition of Theorem \ref{main:thm:1}(1) is now a direct consequence of the definition above.  

\begin{definition}
Let $\F=\{f_1, \ldots, f_n\}$ and $\widetilde{\F}=\{\f_1, \ldots, \f_n\}$ be two systems with  $\sigma(f_i)=\sigma(\f_i)$, for all $1\leq i\leq n$.
Let $K$ be $\F$-invariant and $\widetilde{K}$ be $\widetilde{\F}$-invariant. We say that $\F$ in $K$ is {\it topologically conjugate} to $\widetilde{\F}$ in $\widetilde{K}$ if there exists an order preserving homeomorphism $\varphi:K\to\widetilde{K}$ such that $\varphi(K\cap \textrm{Dom}(f_g))=\widetilde{K}\cap \textrm{Dom}(\f_g)$ for all $g\in G$, and $\varphi\circ f_g(x)=\f_g\circ\varphi(x)$ for all $x\in K\cap \textrm{Dom}(f_g)$.
\end{definition}

Notice that if $K\cap \textrm{Dom}(f_i)$ has at least two points, then the equality $\varphi\circ f_i=\f_i\circ \varphi$ implies that $\sigma(f_i)=\sigma(\f_i)$.

One can also consider the concept of topological semiconjugacy, for this purpose we need to allow degenerate systems.

\begin{definition}
Let $\F=\{f_1, \ldots, f_n\}$ and $\widetilde{\F}=\{\f_1, \ldots, \f_m\}$ be two systems, such that $\widetilde{\F}$ may be degenerate or not. Let $K$ be $\F$-invariant and $\widetilde{K}$ be $\widetilde{\F}$-invariant. We say that $\F$ in $K$ is {\it topologically semiconjugate} to $\widetilde{\F}$ in $\widetilde{K}$, if there exist a map $\tau:\{1,\ldots,n\} \to \{1,\ldots,m\}$ and an order preserving continuous surjective map $\varphi:K\to\widetilde{K}$ such that:
\[ \varphi(K\cap \textrm{Dom}(f_i)) \subset \widetilde{K}\cap \textrm{Dom}(\f_{\tau(i)}),\]
 and $\varphi\circ f_i(x)=\f_{\tau(i)}\circ\varphi(x)$ for all $x\in K\cap \textrm{Dom}(f_i)$ and $1\leq i\leq n$. The map $\varphi$ is called a topological semiconjugacy.
\end{definition}

The map $\tau:\{1,\ldots,n\}\to \{1,\ldots,m\}$ induces a group homomorphism $\tau_*:G\to \widetilde{G}$, where $G$, $\widetilde{G}$ are the groups associated to $\F$, $\widetilde{\F}$ respectively. With this notation we have the following strengthening of the semiconjugacy condition.

\begin{lem}\label{lem:semiG}
Let $\varphi:K\to \widetilde{K}$, $\tau:\{1,\ldots,n\}\to \{1,\ldots,m\}$ give a semiconjugancy between $\F=\{f_1, \ldots, f_n\}$ in $K$ and $\widetilde{\F}=\{\f_1, \ldots, \f_m\}$ in $\widetilde{K}$. Then, for all $g\in G$ we have
\[ \varphi(K\cap {\rm Dom}(f_g)) \subset \widetilde{K}\cap {\rm Dom}(\f_{\tau_*(g)}),\]
 and $\varphi\circ f_g(x)=\f_{\tau_*(g)}\circ\varphi(x)$ for all $x\in K\cap {\rm Dom}(f_g)$.
\end{lem}

\begin{proof}
Notice that if $y\in K\cap \textrm{Dom}(f_i^{-1})=K\cap  \textrm{Im}(f_i)$ then 
$$\varphi(y)\in\varphi(K\cap \textrm{Im}(f_i))\subset \widetilde{K}\cap \f_{\tau(i)}(\varphi(\textrm{Dom}(f_i)\cap K))\subset \widetilde{K}\cap \textrm{Im}(\f_{\tau(i)})=\widetilde{K}\cap \textrm{Dom}(\f_{\tau(i)}^{-1}),$$
we also have
$$\f_{\tau(i)}^{-1}\circ \varphi(y)=\f_{\tau(i)}^{-1}\circ \varphi\circ f_i(f_i^{-1}(y))=\f_{\tau(i)}^{-1}\circ\f_{\tau(i)}\circ\varphi(f_i^{-1}(y))=\varphi\circ f_i^{-1}(y).$$
This proves the desired properties for $|g|=1$. For the rest of the proof we will proceed by induction on $|g|$.

Assume the case $|g|=k$ is true and we will show $|g|=k+1$. If $g=a_it$ with $|t|=k$, then we have that $x\in K\cap\textrm{Dom}(f_g)$ implies $f_{a_i}(x)\in K\cap\textrm{Dom}(f_t)$ and by the induction hypothesis
\[\varphi(f_{a_i}(x))\in\widetilde{K}\cap\textrm{Dom}(\f_{\tau_*(t)}), \text{ and } \varphi(f_t(f_{a_i}(x)))=\f_{\tau_*(t)}(\varphi(f_{a_i}(x))).\]
This two properties transform into
\[\f_{\tau_*(a_i)}(\varphi(x))\in\widetilde{K}\cap \textrm{Dom}(\f_{\tau_*(t)}), \text{ and } \varphi(f_t(f_{a_i}(x)))=\f_{\tau_*(t)}(\f_{\tau_*(a_i)}(\varphi(x))),\]
which are the desired properties for $g=a_it$. The case in which $g$ is of the form $a_i^{-1}t$ is analogous. This finishes the proof of the lemma.
\end{proof}



The proposition below ensures that the kneading sequences 
are preserved by topological conjugacy. 

\begin{prop}\label{prop:top}
Let $\F=\{f_1, \ldots, f_n\}$ and $\widetilde{\F}=\{\f_1, \ldots, \f_n\}$ be two systems 
with ${\rm Dom}(f_i)={\rm Dom}(\f_i)$ 
and $\sigma(f_i)=\sigma(\f_i)$, for all $1\leq i\leq n$.
If $\F$ and $\widetilde{\F}$ are topologically conjugate at invariant sets containing the critical $G$-orbits, then they have the same kneading sequences.
\end{prop}

\begin{proof}
We have to show that $S_c=\widetilde{S}_c$ and $T_{c}\rr_{S_c}=\widetilde{T}_c\rr_{\widetilde{S}_c}$
for all turning point $c$. Let $\varphi:K\to\widetilde{K}$ be an order preserving homeomorphism conjugating $\F$ in $K$ and $\widetilde{\F}$ in $\widetilde{K}$, where $K,\widetilde{K}$ are invariant sets containing $\cc_G,\widetilde{\cc}_G$, respectively. So $\varphi(c_i^\ell)=\varphi(\inf(K\cap \textrm{Dom}(f_i)))=
\inf\varphi(K\cap \textrm{Dom}(f_i))= \inf(\widetilde{K}\cap \textrm{Dom}(\widetilde{f}_i))=c_i^\ell$, and analogously $\varphi(c_i^r)=c_i^r$. Thus $\varphi(c)=c$ for all turning point, hence $S_c=\widetilde{S}_c$ (indeed $G_c=\widetilde{G}_c$).  
Since $\varphi$ is order preserving and $\varphi(f_g(c))=\f_g(c)$ for all $c$ turning point and $g\in G_c$ we have
\[
\lambda(f_g(c),c')=\lambda(\f_g(c),c'),\ \text{ for all } c' \text{ turning point}.
\]
Therefore $T_{c}\rr_{G_c}=\widetilde{T}_c\rr_{\widetilde{G}_c}$.
In particular $T_{c}\rr_{S_c}=\widetilde{T}_c\rr_{\widetilde{S}_c}$. 
\end{proof}

In the next propositions we will study how the kneading sequences influence the combinatorics of the orbits of turning points.
Given $x,\, y \in \R$, we denote by $[x,y]$ the closed interval of real numbers between $x$ and $y$ (note that $[x,y]=[y,x]$).

\begin{prop}\label{prop:fp}
Let $\F=\{f_1, \ldots, f_n\}$ and $\widetilde{\F}=\{\widetilde{f}_1, \ldots, \widetilde{f}_n\}$ be two systems with ${\rm Dom}(f_i)={\rm Dom}(\widetilde{f}_i)$ and $\sigma(f_i)=\sigma(\widetilde{f}_i)$, for all $1\leq i\leq n$. If for all turning points $c$ we have $S_c=\widetilde{S}_c$ 
and $T_{c}\rr_{S_c}=\widetilde{T}_c\rr_{\widetilde{S}_c}$, then $G_c\cap S^{-1}=\widetilde{G}_c\cap S^{-1}$ and
$$
\lambda(f_g(c),f_h(c'))=\lambda(\widetilde{f}_g(c),\widetilde{f}_h(c')),
$$
for all turning points $c, c'$ and $h\in S_{c'}=\widetilde{S}_{c'}$, $g\in G_c\cap S^{-1}=\widetilde{G}_c\cap S^{-1}$. In particular $T_{c}\rr_{G_c\cap S^{-1}}=\widetilde{T}_c\rr_{\widetilde{G}_c\cap S^{-1}}$ (we get this putting $h=e$).
\end{prop}

\begin{proof}
Let $g\in G_c\cap S^{-1}$, we will prove that $g\in \widetilde{G}_c\cap S^{-1}$ and that for any turning point $c'$ and $h\in S_{c'}=\widetilde{S}_{c'}$ we have
$$
\lambda(f_g(c),f_h(c'))=\lambda(\widetilde{f}_g(c),\widetilde{f}_h(c')).
$$
We proceed by induction on $|g|$, the base case, $|g|=0$, corresponds to the hypothesis in the proposition. Write $g=ta_i^{-1}$ with $|t|=|g|-1$, since $g\in G_c\cap S^{-1}$ then $t\in G_c\cap S^{-1}$ and by the induction hypothesis $t\in \widetilde{G}_c\cap S^{-1}$. Moreover
\[\lambda(f_t(c),f_{a_i}(c_i^\ell))=\lambda(\widetilde{f}_t(c),\widetilde{f}_{a_i}(c_i^\ell)),\,\,\lambda(f_t(c),f_{a_i}(c_i^r))=\lambda(\widetilde{f}_t(c),\widetilde{f}_{a_i}(c_i^r)),\]
therefore $f_t(c)$ is inside the interval $[f_{a_i}(c_i^\ell),f_{a_i}(c_i^r)]$ if and only if $\widetilde{f}_t(c)$ is inside $[\widetilde{f}_{a_i}(c_i^\ell),\widetilde{f}_{a_i}(c_i^r)]$. Since $g=ta_i^{-1}\in G_c\cap S^{-1}$ then $g\in \widetilde{G}_c\cap S^{-1}$. 

Now if $h a_i \notin S_{c'}=\widetilde{S}_{c'}$ then $f_h(c')$ is outside of the interval $I_i$, the same would happen with $\widetilde{f}_h(c')$. Moreover, since $D(f_h(c'))=D(\widetilde{f}_h(c'))$ then $f_h(c')$ is to the left or right of $I_i$ if and only if $\widetilde{f}_h(c')$ is. Given that $f_g(c),\widetilde{f}_g(c) \in I_i$ we get that $\lambda(f_g(c),f_h(c'))=\lambda(\widetilde{f}_g(c),\widetilde{f}_h(c'))$.

On the other hand, if $h a_i \in S_{c'}=\widetilde{S}_{c'}$ then by the induction hypothesis we have
\begin{align*}
\lambda(f_g(c),f_h(c'))&=\sigma(f_{a_i}) \lambda(f_t(c),f_{h a_i}(c'))\\
&=\sigma(\widetilde{f}_{a_i}) \lambda(\widetilde{f}_{t}(c),\widetilde{f}_{h a_i}(c'))\\
&=\lambda(\widetilde{f}_g(c),\widetilde{f}_h(c')).
\end{align*}
We have proven by induction that $G_c\cap S^{-1} \subset \widetilde{G}_c\cap S^{-1}$ and
\[\lambda(f_g(c),f_h(c'))=\lambda(\widetilde{f}_g(c),\widetilde{f}_h(c')),\]
by symmetry $G_c\cap S^{-1} = \widetilde{G}_c\cap S^{-1}$ as we wanted.
\end{proof}

With the right hypothesis we can prove stronger results than Proposition \ref{prop:fp}.

\medskip
\noindent
{\sc Separation property:} Given two points $x,y$ we say that an element $g\in G$ separates them if one of the following holds:
\begin{enumerate}[{\rm (i)}]
\item $x\in \textrm{Dom}(f_g)$ and $y\notin \textrm{Dom}(f_g)$.
\item $y\in \textrm{Dom}(f_g)$ and $x\notin \textrm{Dom}(f_g)$.
\item $x,y \in \textrm{Dom}(f_g)$ and $D(f_g(x))\neq D(f_g(y))$.
\end{enumerate}


\begin{rem}
If $t\in S$ separates $x,\, y$,
then we may assume that (iii) occurs. Indeed, let $t\in S$ and
suppose that $x\in\textrm{Dom}(f_t)$ and $y\notin \textrm{Dom}(f_t)$
(the case $y\in\textrm{Dom}(f_t)$ and $x\notin\textrm{Dom}(f_t)$
is treated similarly). Consider a sequence $f_{t_0}, \ldots, f_{t_j}$ 
so that $t_0$ is the identity of $S$, $t_j=t$ and $t_{m}=t_{m-1} a_{i(m)}$ for 
any $1\leq m\leq j$ and some generating element $a_{i(m)}$. 
Note that there exists $1\leq m\leq j$ such that $y\in\textrm{Dom}
(f_{t_{m-1}})$ and $y\notin\textrm{Dom}(f_{t_{m}})$.
Hence $f_{t_{m-1}}(y)$ is either to the right or left of $I_{i(m)}$.
Since $f_{t_{m-1}}(x)\in I_{i(m)}$, it follows that 
$D(f_{t_{m-1}}(x))\neq D(f_{t_{m-1}}(y))$. So just
take $t=t_{m-1}$.
\end{rem}

\begin{definition}
Given a system $\F=\{f_1, \ldots, f_n\}$, we will say that:
\begin{enumerate}[$\circ$]
\item $\F$ {\it separates points in the past}, if for any two different generators $a,\, b$  and any two different points $x\in f_{a^{-1}}(\cc_{S^{-1}})$, $y\in f_{b^{-1}}(\cc_{S^{-1}})$, $x,\, y$ can be separated by an element in $S^{-1}$.
\item $\F$ {\it separates the critical orbit in the future}, if any two different points $x,\, y \in \cc_G$ can be separated by an element in $S$.
\item $\F$ verifies {\it the separability hypothesis} if it separates points in the past and separates the critical orbit in the future.
\end{enumerate}
\end{definition}

For examples of systems that satisfy the separability hypothesis see Section \ref{appexamples}. 

\begin{prop}\label{prop:pp}
Assume the same setting and hypothesis as Proposition \ref{prop:fp}. In addition, suppose that $\F$ separates points in the past. Then we have
$$
\lambda(f_g(c),f_t(c'))\neq 0 \Rightarrow \lambda(f_g(c),f_t(c'))=\lambda(\widetilde{f}_g(c),\widetilde{f}_t(c')),
$$
for all turning points $c, c'$ and $t\in G_{c'}\cap S^{-1}=
\widetilde{G}_{c'}\cap S^{-1}$, $g\in G_c\cap S^{-1}=\widetilde{G}_c\cap S^{-1}$. Moreover, if we assume that $\widetilde{\F}$ also separates points in the past then
$$
\lambda(f_g(c),f_t(c'))=\lambda(\widetilde{f}_g(c),\widetilde{f}_t(c')),
$$
for all turning points $c, c'$ and $t\in G_{c'}\cap S^{-1}=
\widetilde{G}_{c'}\cap S^{-1}$, $g\in G_c\cap S^{-1}=\widetilde{G}_c\cap S^{-1}$.
\end{prop}
\begin{proof}
Let $t\in G_{c'} \cap S^{-1}=\widetilde{G}_{c'}\cap S^{-1}$ and $g\in G_c\cap 
S^{-1}=\widetilde{G}_c\cap S^{-1}$. We will proceed by induction 
in $|g|+|t|$, the base case $|g|+|t|=0$ is trivial. Notice 
that if $|g|=0$ or $|t|=0$ then the desired equality follows 
from Proposition \ref{prop:fp}. We can assume $|g|\geq 1$, 
$|t|\geq 1$ and write $g=h_1 a_i^{-1}$, $t=h_2 a_j^{-1}$ 
with $|h_1|=|g|-1$, $|h_2|=|t|-1$. If $a_i=a_j$ and $\lambda(f_g(c),f_t(c'))\neq 0$, using the induction hypothesis we have
\begin{align*}
\lambda(f_g(c),f_t(c'))&= \sigma(f_{a_i}) \lambda(f_{h_1}(c),f_{h_2}(c'))\\
&=\sigma(\widetilde{f}_{a_i}) \lambda(\widetilde{f}_{h_1}(c),\widetilde{f}_{h_2}(c'))\\
&=\lambda(\widetilde{f}_g(c),\widetilde{f}_t(c')).
\end{align*}
If $a_i\neq a_j$ and $\lambda(f_g(c),f_t(c'))\neq 0$ then $f_g(c)\in f_{a_i^{-1}}(\cc_{S^{-1}})$, $f_t(c')\in f_{a_j^{-1}}(\cc_{S^{-1}})$ and from the hypothesis we get that they can be separated by an element in $S^{-1}$. We have three cases:
\begin{enumerate}[$\circ$]
\item There is $J\in S^{-1}$ and a generator $a_l$ such that $f_{g}(c) \in \text{Dom}(f_{Ja_l^{-1}})$ and $f_{t}(c')\in \text{Dom}(f_J)\setminus \text{Dom} (f_{Ja_l^{-1}})$. In this case $f_{gJ}(c)$ is in the interval $[f_{a_l}(c_l^\ell),f_{a_l}(c_l^r)]$ and $f_{tJ}(c')$ is outside of it. Now, Proposition \ref{prop:fp} implies that 
$\lambda(f_{gJ}(c),f_{a_l}(c_l^{\ell/r}))=\lambda(\widetilde{f}_{gJ}(c),\widetilde{f}_{a_l}(c_l^{\ell/r}))$ and $\lambda(f_{tJ}(c'),f_{a_l}(c_l^{\ell/r}))=\lambda(\widetilde{f}_{tJ}(c'),\widetilde{f}_{a_l}(c_l^{\ell/r}))$, and so
\begin{align*}
\lambda(f_g(c),f_t(c'))&= \sigma(f_J) \lambda(f_{gJ}(c),f_{tJ}(c'))
=\sigma(f_J) \lambda(f_{a_l}(c_l^r),f_{tJ}(c'))\\
&=\sigma(\widetilde{f}_J)\lambda(\widetilde{f}_{a_l}(c_l^r),\widetilde{f}_{tJ}(c))=\sigma(\widetilde{f}_J)\lambda(\widetilde{f}_{gJ}(c),\widetilde{f}_{tJ}(c'))\\
&=\lambda(\widetilde{f}_g(c),\widetilde{f}_t(c')).
\end{align*}
\item There is $J\in S^{-1}$ and a generator $a_l$ such that $f_{t}(c) \in \text{Dom}(f_{Ja_l^{-1}})$ and $f_{g}(c')\in \text{Dom}(f_J)\setminus \text{Dom} (f_{Ja_l^{-1}})$. This case is similar to the previous one.
\item There is $J\in S^{-1}$ such that $f_{g}(c),\, f_{t}(c') \in \text{Dom}(f_{J})$ and $D(f_{tJ}(c'))\neq D(f_{gJ}(c))$. By proposition \ref{prop:fp} we have $D(f_{tJ}(c'))=D(\widetilde{f}_{tJ}(c'))$ and $D(f_{gJ}(c))=D(\widetilde{f}_{gJ}(c))$. In this case, Lemma \ref{lem:lamb} implies that $\lambda(f_{tJ}(c'),f_{gJ}(c))$ only depends on $D(f_{tJ}(c'))$, $D(f_{gJ}(c))$, and hence
\begin{align*}
\lambda(f_{t}(c'),f_{g}(c))&=\sigma(f_J) \lambda(f_{tJ}(c'),f_{gJ}(c))\\
&=\sigma(\widetilde{f}_J)\lambda(\widetilde{f}_{tJ}(c'),\widetilde{f}_{gJ}(c))=\lambda(\widetilde{f}_{t}(c'),\widetilde{f}_{g}(c)).
\end{align*}
\end{enumerate}

To finish, assume both systems separate points in the past.
Then $\lambda(f_g(c),f_t(c'))\neq 0$ or $\lambda(\widetilde{f}_g(c),\widetilde{f}_t(c'))\neq 0$  would imply $\lambda(f_g(c),f_t(c'))=\lambda(\widetilde{f}_g(c),\widetilde{f}_t(c'))$, otherwise we have $\lambda(f_g(c),f_t(c'))=\lambda(\widetilde{f}_g(c),\widetilde{f}_t(c'))= 0$. Thus in both cases we get the desired equality.
\end{proof}

In the particular case in which the intervals intersect only at border points, the additional assumption in the previous lemma is verified. In fact, in this case one has that if $D(f_{a_i^{-1}}(x))=
D(f_{a_j^{-1}}(y))$, for $i\neq j$, then $f_{a_i^{-1}}(x)=
f_{a_j^{-1}}(y)$ and they are equal to the turning point in 
the intersection between $I_i$ and $I_j$. Thus, the identity element $e$ already separates two different points $f_{a_i^{-1}}(x)$, $f_{a_j^{-1}}(y)$. 

\begin{cor}\label{cor:iob}
Assume the same setting and hypothesis as Proposition \ref{prop:fp}. 
In addition, suppose that for any two intervals in $I_1, \ldots, I_n$, 
the intersection is either empty or a turning point. Then we have
$$
\lambda(f_g(c),f_t(c'))=\lambda(\widetilde{f}_g(c),\widetilde{f}_t(c')),
$$
for all turning points $c, c'$ and $t\in G_{c'}\cap S^{-1}=
\widetilde{G}_{c'}\cap S^{-1}$, $g\in G_c\cap S^{-1}=\widetilde{G}_c\cap S^{-1}$.
\end{cor}

Proposition \ref{prop:fp} showed that if two systems have the same kneading sequences, then the relative position between pairs of points, one in the $S$--orbit and one in the $S^{-1}$--orbit, is the same for the two systems. Our objective will be to obtain this property for pairs of points in the whole $G$--orbit, but this is not true in general. Example \ref{exa:nsp} gives two systems with the same kneading sequences but such that pairs of points, both in the $S^{-1}$--orbit, do not have the same relative position for the two systems. This is why we needed to introduce the separability hypothesis to prove Proposition \ref{prop:pp}.

Example \ref{exp:ff} shows that even if both systems separate points in the past, then one can not guaranteed that pairs of points in the $S$--orbit have the same relative position for the two systems. We will be able to prove same relative position for the whole $G$--orbit assuming both system satisfy the separability hypothesis, this is the content of Theorem \ref{thm:ft}.


\begin{example}\label{exa:nsp}
Let $\F=\{f_1,f_2\}$, $\widetilde{\F}=\{\widetilde{f}_1,\widetilde{f}_2\}$ be two systems 
defined on the intervals $I_1,I_2$. 
These two intervals are equal $I_1=I_2=[-1,1]$ and the functions are defined as 
\begin{align*}
&f_1(x)=2x, \, f_2(x)=3x, \\
&\widetilde{f}_1(x)=2x, \, \widetilde{f}_2(x)=3x+\tfrac{3}{4}.
\end{align*}
Since $[-1,1]\subset{\rm Dom}(f_g)$, $[-1,1]\subset{\rm Dom}(\f_g)$ and $f_g(\pm1),\widetilde{f}_g(\pm1)\in(-1,1)$ for all $g\in S^{-1}$, we have that both $\F,\widetilde{\F}$ do not separate points in the past.  
It is easy to see that $S_{-1}=\widetilde{S}_{-1}$, $S_{1}=\widetilde{S}_{1}$ 
and $D(f_g(-1))=D(\widetilde{f}_g(-1))$, $D(f_h(1))=D(\widetilde{f}_h(1))$ 
for all $g\in S_{-1}=\widetilde{S}_{-1}$, $h\in S_{1}=\widetilde{S}_{1}$. 
However $f_1^{-1}(-1)< f_2^{-1}(-1)$ and $\widetilde{f}_1^{-1}(-1)> \widetilde{f}_2^{-1}(-1)$. 
\end{example}

\begin{example}\label{exp:ff}
Let $\F=\{f_1,f_2\}$, $\widetilde{\F}=\{\widetilde{f}_1,\widetilde{f}_2\}$ be two systems 
defined on the intervals $I_1=[-1,0]$ and $I_2=[0,1]$. 
The functions are defined as 
\begin{align*}
&f_1(x)=\tfrac{1}{2}(x+1), \, f_2(x)=\tfrac{1}{4}(3x+1), \\
&\widetilde{f}_1(x)=\tfrac{1}{2}(x+1), \,  \widetilde{f}_2(x)=\tfrac{1}{4}(x+3).
\end{align*}
Again it is easy to see that $S_{-1}=\widetilde{S}_{-1}$, $S_{0}=\widetilde{S}_{0}$, $S_{1}=\widetilde{S}_{1}$ 
and $D(f_g(-1))=D(\widetilde{f}_g(-1))$, $D(f_s(0))=D(\widetilde{f}_s(0))$, $D(f_h(1))=D(\widetilde{f}_h(1))$ 
for all $g\in S_{-1}=\widetilde{S}_{-1}$, $s\in S_{0}=\widetilde{S}_{0}$,  $h\in S_{1}=\widetilde{S}_{1}$. 
Nevertheless $f_1(0)>f_2(0)$ and $\widetilde{f}_1(0)< \widetilde{f}_2(0)$.
\end{example}

\begin{thm}\label{thm:ft}
Let $\F=\{f_1, \ldots, f_n\}$ and $\widetilde{\F}=\{\widetilde{f}_1, \ldots, \widetilde{f}_n\}$ 
be two systems with ${\rm Dom}(f_i)={\rm Dom}(\widetilde{f}_i)$ and 
$\sigma(f_i)=\sigma(\widetilde{f}_i)$, for all $1\leq i\leq n$. 
Assume that both systems satisfy the separability hypothesis. If for all turning points $c$ we have $S_c=\widetilde{S}_c$ 
and $T_{c}\rr_{S_c}=\widetilde{T}_c\rr_{\widetilde{S}_c}$, then $G_c=\widetilde{G}_c$ and
\begin{equation}\label{eq:ordem}
\lambda(f_g(c),f_h(c'))=\lambda(\widetilde{f}_g(c),\widetilde{f}_h(c')),
\end{equation}
for all turning points $c, c'$ and $h\in G_{c'}=\widetilde{G}_{c'}$, $g\in G_c=\widetilde{G}_c$. In particular $T_{c}=\widetilde{T}_c$ (we get this putting $h=e$).
\end{thm}
\begin{proof}
Every $g\in G\setminus\{e\}$ can be written as $g=B_1\ldots B_l$, where every $B_i$ is in $S\setminus \{e\}$ or $S^{-1}\setminus \{e\}$ and two consecutive $B_i's$ belong to different sets. For a fixed value of $l$ we define:
\begin{align*}
Y_l^+&=\{g=B_1\ldots B_l:\, B_l \in S\},\\
Y_l^-&=\{g=B_1\ldots B_l:\, B_l \in S^{-1}\}.
\end{align*}
We will prove the following statement: for any turning points $c,c'$ we have $Y_l^+\cap G_c=Y_l^+\cap \widetilde{G}_c$, $Y_l^-\cap G_c=Y_l^-\cap \widetilde{G}_c$ and
\begin{align*}
\lambda(f_g(c),f_t(c'))&=\lambda(\widetilde{f}_g(c),\widetilde{f}_t(c')),\\
\lambda(f_h(c),f_s(c'))&=\lambda(\widetilde{f}_h(c),\widetilde{f}_s(c')),
\end{align*}
for all $g \in Y_l^+\cap G_c$, $h\in Y_l^-\cap G_c$, $t\in S^{-1}\cap G_{c'}$, $s\in S_{c'}$.

The case $l=1$ follows from Proposition \ref{prop:fp}. We will proceed by induction on $l$, and for a fixed value of $l$ we will do induction on $|B_l|$. Let $g=B_1\ldots B_l\in Y_l^+\cap G_c$, we can write $B_l=A a_i$ with $|A|=|B_l|-1$ and $A\in S$. If $A=e$ then $B_1\ldots B_{l-1} \in Y_{l-1}^-\cap G_c$ and by the induction hypothesis $B_1\ldots B_{l-1} \in Y_{l-1}^-\cap \widetilde{G}_c$. Otherwise $B_1\ldots B_{l-1}A \in Y_l^+\cap G_c$ and by the induction hypothesis $B_1\ldots B_{l-1}A \in Y_l^+\cap \widetilde{G}_c$. In both cases we have $D(f_{B_1\ldots B_{l-1}A}(c))=D(\widetilde{f}_{B_1\ldots B_{l-1}A}(c))$
 and then $\widetilde{f}_{B_1\ldots B_{l-1}A}(c)\in I_i$ which implies $g\in Y_l^+\cap \widetilde{G}_c$. This proves $Y_l^+\cap G_c \subset Y_l^+\cap \widetilde{G}_c$, analogously one proves the other way around, thus $Y_l^+\cap G_c = Y_l^+\cap \widetilde{G}_c$.

On the other hand, given $t\in S^{-1}\cap G_{c'}$ we have two cases:

\noindent
{\sc Case 1:} If $f_t(c') \in \textrm{Dom}(f_{a_i^{-1}})$ then 
$\widetilde{f}_t(c') \in \textrm{Dom}(\widetilde{f}_{a_i^{-1}})$ 
(thanks to Proposition \ref{prop:fp}), and we have two options: 
\begin{enumerate}[$\circ$]
\item $|A|\geq 1$: Using the induction hypothesis we get
\begin{align*}
\lambda(f_g(c),f_t(c'))&=\sigma(f_{a_i}) \lambda(f_{B_1\ldots B_{l-1}A}(c),f_{ta_i^{-1}}(c')),\\
&=\sigma(\widetilde{f}_{a_i}) \lambda(\widetilde{f}_{B_1\ldots B_{l-1}A}(c),\widetilde{f}_{ta_i^{-1}}(c')),\\
&=\lambda(\widetilde{f}_g(c),\widetilde{f}_t(c')).
\end{align*}

\item $A=e$: If $l=2$, Proposition \ref{prop:pp} implies $\lambda(f_{B_1}(c),f_{ta_i^{-1}}(c'))=\lambda(\widetilde{f}_{B_1}(c),\widetilde{f}_{ta_i^{-1}}(c'))$ and thus $\lambda(f_g(c),f_t(c'))=\lambda(\widetilde{f}_g(c),\widetilde{f}_t(c'))$. Assuming now that $l>2$ we have $f_{B_1\ldots B_{l-1}}(c)\in \text{Dom}(f_{B_{l-1}^{-1}})$,  
and the same for the other system. We have two alternatives:
\begin{enumerate}[{\rm (1)}]
\item If $f_{ta_i^{-1}}(c')\in \text{Dom}(f_{B_{l-1}^{-1}})$ then $\widetilde{f}_{ta_i^{-1}}(c')\in \text{Dom}(\widetilde{f}_{B_{l-1}^{-1}})$ (induction hypothesis). Using the separation condition we know that if $f_{ta_i^{-1}B_{l-1}^{-1}}(c')\neq f_{B_1\ldots B_{l-2}}(c)$ then there is $J\in S$ such that $D(f_{ta_i^{-1}B_{l-1}^{-1}J}(c'))\neq D(f_{B_1\ldots B_{l-2}J}(c))$. By the induction hypothesis (remember that $l>2$) we have that $D(f_{ta_i^{-1}B_{l-1}^{-1}J}(c'))=D(\widetilde{f}_{ta_i^{-1}B_{l-1}^{-1}J}(c'))$ and $D(f_{B_1\ldots B_{l-2}J}(c))=D(\widetilde{f}_{B_1\ldots B_{l-2}J}(c))$, therefore
\begin{equation}\label{eq:se}
\lambda(f_{B_1\ldots B_{l-2}J}(c),f_{ta_i^{-1}B_{l-1}^{-1}J}(c'))=\lambda(\widetilde{f}_{B_1\ldots B_{l-2}J}(c),\widetilde{f}_{ta_i^{-1}B_{l-1}^{-1}J}(c')).
\end{equation}
A similar argument can be used to get the same equality if $\widetilde{f}_{ta_i^{-1}B_{l-1}^{-1}}(c')\neq \widetilde{f}_{B_1\ldots B_{l-2}}(c)$, otherwise we have
\[\lambda(f_{B_1\ldots B_{l-2}}(c),f_{ta_i^{-1}B_{l-1}^{-1}}(c'))=\lambda(\widetilde{f}_{B_1\ldots B_{l-2}}(c),\widetilde{f}_{ta_i^{-1}B_{l-1}^{-1}}(c'))=0.\]
Either way, we get equation (\ref{eq:se}) for some $J\in S$. We can conclude
\begin{align*}
\lambda(f_g(c),f_t(c'))&=\sigma(f_{a_i}) \sigma(f_{B_{l-1}}) \sigma(f_J) \lambda(f_{B_1\ldots B_{l-2}J}(c),f_{ta_i^{-1}B_{l-1}^{-1}J}(c'))\\
&=\sigma(\widetilde{f}_{a_i}) \sigma(\widetilde{f}_{B_{l-1}}) \sigma(\widetilde{f}_J) \lambda(\widetilde{f}_{B_1\ldots B_{l-2}J}(c),\widetilde{f}_{ta_i^{-1}B_{l-1}^{-1}J}(c'))\\
&=\lambda(\widetilde{f}_g(c),\widetilde{f}_t(c')).
\end{align*}
\item If $f_{ta_i^{-1}}(c')\notin \text{Dom}(f_{B_{l-1}^{-1}})$ then $B_{l-1}^{-1}$ separates $f_{ta_i^{-1}}(c')$ and $f_{B_1\ldots B_{l-1}}(c)$. Then there is $J\in S$ such that $B_{l-1}^{-1}$ starts with $J$ and $D(f_{ta_i^{-1}J}(c'))\neq D(f_{B_1\ldots B_{l-1}J}(c))$. By the induction hypothesis (remember $l>2$ and $B_{l-1}J\in S^{-1}$) we have $D(f_{B_1\ldots B_{l-1}J}(c))=D(\widetilde{f}_{B_1\ldots B_{l-1}J}(c))$ and  $D(f_{ta_i^{-1}J}(c'))=D(\widetilde{f}_{ta_i^{-1}J}(c'))$, therefore
\[\lambda(f_{ta_i^{-1}J}(c'),f_{B_1\ldots B_{l-1}J}(c))=\lambda(\widetilde{f}_{ta_i^{-1}J}(c'),\widetilde{f}_{B_1\ldots B_{l-1}J}(c)).\]
We can conclude
\begin{align*}
\lambda(f_g(c),f_t(c'))&=\sigma(f_{a_i}) \sigma(f_J) \lambda(f_{B_1\ldots B_{l-1}J}(c),f_{ta_i^{-1}J}(c'))\\
&=\sigma(\widetilde{f}_{a_i}) \sigma(\widetilde{f}_J) \lambda(\widetilde{f}_{B_1\ldots B_{l-1}J}(c),\widetilde{f}_{ta_i^{-1}J}(c'))\\
&=\lambda(\widetilde{f}_g(c),\widetilde{f}_t(c')).
\end{align*}
\end{enumerate}
\end{enumerate}

\noindent
{\sc Case 2:} If $f_t(c') \notin\textrm{Dom}(f_{a_i^{-1}})$ then 
$f_t(c')$ is outside $[f_{a_i}(c_i^\ell),f_{a_i}(c_i^r)]$ but $f_g(c)$ is inside of it. Therefore, using the induction hypothesis, we get
\begin{align*}
\lambda(f_g(c), f_t(c'))&=\lambda(f_{a_i}(c_i^r),f_t(c'))\\
&=\lambda(\widetilde{f}_{a_i}(c_i^r),\widetilde{f}_t(c'))\\
&=\lambda(\widetilde{f}_g(c), \widetilde{f}_t(c')).
\end{align*}

It remains to consider $g=B_1\ldots B_l\in Y_l^-\cap G_c$. This is similar to the treatment for $Y^+_l$, for completeness we present the proof. We can write $B_l=A a_i^{-1}$ with $|A|=|B_l|-1$ and $A\in S^{-1}$. Using the induction hypothesis and the separation condition we can prove that
\begin{align*}
\lambda(f_{B_1\ldots B_{l-1}A}(c),f_{a_i}(c_i^r))&=\lambda(\widetilde{f}_{B_1\ldots B_{l-1}A}(c),\widetilde{f}_{a_i}(c_i^r)),\\
\lambda(f_{B_1\ldots B_{l-1}A}(c),f_{a_i}(c_i^\ell))&=\lambda(\widetilde{f}_{B_1\ldots B_{l-1}A}(c),\widetilde{f}_{a_i}(c_i^\ell)).
\end{align*}
The equalities are straight forward from the induction hypothesis in the case that $|A|\geq 1$, if $A=e$ then to prove the equalities we also need to consider an element of $S$ separating $f_{B_1\ldots B_{l-1}}(c)$ and $f_{a_i}(c_i^r)$ (or $f_{a_i}(c_i^\ell)$ for the other equality), this can easily be done. Using these equalities we can conclude that $\widetilde{f}_{B_1\ldots B_{l-1}A}(c)\in {\rm Dom}(\f_{a_i^{-1}})$ and this proves $g\in Y_l^-\cap \widetilde{G}_c$. We have proven $Y_l^-\cap G_c \subset Y_l^-\cap \widetilde{G}_c$,
analogously one proves the other way around, thus $Y_l^-\cap G_c = Y_l^-\cap
\widetilde{G}_c$.

On the other hand, given $s\in S_{c'}$ we have two cases: 

\noindent
{\sc Case 1:} If $f_s(c') \in \textrm{Dom}(f_{a_i})$ then 
$\widetilde{f}_s(c') \in \textrm{Dom}(\widetilde{f}_{a_i})$, and we have two options: 
\begin{enumerate}[$\circ$]
\item $|A|\geq 1$: Using the induction hypothesis we get
\begin{align*}
\lambda(f_g(c),f_s(c'))&=\sigma(f_{a_i}) \lambda(f_{B_1\ldots B_{l-1}A}(c),f_{s a_i}(c')),\\
&=\sigma(\widetilde{f}_{a_i}) \lambda(\widetilde{f}_{B_1\ldots B_{l-1}A}(c),\widetilde{f}_{s a_i}(c')),\\
&=\lambda(\widetilde{f}_g(c),\widetilde{f}_s(c')).
\end{align*}

\item $A=e$: If $f_{B_1\ldots B_{l-1}}(c)\neq  
f_{sa_i}(c')$, we use the separation hypothesis to find an element of $J\in S$ such that $f_{B_1\ldots B_{l-1}}(c),  
f_{sa_i}(c') \in \textrm{Dom}(f_J)$ and $D(f_{B_1\ldots B_{l-1}J}(c))\neq D(f_{sa_iJ}(c'))$. In this case, $\lambda(f_{B_1\ldots B_{l-1}J}
(c),f_{sa_iJ}(c'))$ is determined by the value of $D(f_{B_1\ldots 
B_{l-1}J}(c)),\, D(f_{sa_iJ}(c'))$. Moreover, the induction 
hypothesis implies that $D(f_{B_1\ldots B_{l-1}J}(c))=D(\widetilde{f}_{B_1\ldots 
B_{l-1}J}(c))$ and $D(f_{sa_iJ}(c'))=D(\widetilde{f}_{sa_iJ}(c'))$, 
which gives \[\lambda(f_{B_1\ldots B_{l-1}J}(c),f_{sa_iJ}(c'))=
\lambda(\widetilde{f}_{B_1\ldots B_{l-1}J}(c),\widetilde{f}_{sa_iJ}(c')).\]
If $\widetilde{f}_{B_1\ldots B_{l-1}}(c)\neq\widetilde{f}_{sa_i}(c')$ we use a similar argument to find the same equality. Otherwise we have
\[\lambda(f_{B_1\ldots B_{l-1}}(c),f_{sa_i}(c'))=
\lambda(\widetilde{f}_{B_1\ldots B_{l-1}}(c),\widetilde{f}_{sa_i}(c'))=0.\]
Either way, one gets
$\lambda(f_{B_1\ldots B_{l-1}J}(c),f_{sa_iJ}(c'))=
\lambda(\widetilde{f}_{B_1\ldots B_{l-1}J}(c),\widetilde{f}_{sa_iJ}(c'))$
for some $J\in S$. We conclude that
\begin{align*}
\lambda(f_g(c),f_s(c'))&=\sigma(f_{a_i})\lambda(f_{B_1\ldots B_{l-1}}(c),f_{sa_i}(c'))\\
&=\sigma(f_{a_i})\sigma(f_J)\lambda(f_{B_1\ldots B_{l-1}J}(c),f_{sa_iJ}(c'))\\
&=\sigma(\widetilde{f}_{a_i})\sigma(\widetilde{f}_J)\lambda(\widetilde{f}_
{B_1\ldots B_{l-1}J}(c),\widetilde{f}_{sa_iJ}(c'))\\
&=\lambda(\widetilde{f}_g(c),\widetilde{f}_s(c')).
\end{align*}
\end{enumerate}

\noindent
{\sc Case 2:} If $f_s(c') \notin\textrm{Dom}(f_{a_i})$ then 
$f_s(c')$ is outside $[c_i^\ell,c_i^r]$ but 
$f_g(c)$ is inside of it. Therefore
\begin{align*}
\lambda(f_g(c), f_t(c'))&=\lambda(c_i^r,f_s(c'))\\
&=\lambda(c_i^r,\widetilde{f}_s(c'))\\
&=\lambda(\widetilde{f}_g(c), \widetilde{f}_s(c')).
\end{align*}

Now that we have proven the statement we can prove the theorem. Notice that thanks to the statement we already have that $G_c=\widetilde{G}_c$ and $D(f_g(c))=D(\widetilde{f}_g(c))$ for any turning point $c$ and  $g\in G_c$. Take another turning point $c'$ and $h\in G_{c'}=\widetilde{G}_{c'}$, if $f_g(c)\neq f_h(c')$ by the separation hypothesis there is $J\in S$ such that $f_g(c),f_h(c')\in {\rm Dom}(f_J)$ and $D(f_{gJ}(c))\neq D(f_{hJ}(c'))$. Using that $D(f_{gJ}(c))=D(\widetilde{f}_{gJ}(c))$, $D(f_{hJ}(c'))=D(\widetilde{f}_{hJ}(c'))$ we conclude that
\begin{align*}
\lambda(f_g(c),f_h(c'))&=\sigma(f_J) \lambda(f_{gJ}(c),f_{hJ}(c'))\\
&=\sigma(\widetilde{f}_J) \lambda(\widetilde{f}_{gJ}(c),\widetilde{f}_{hJ}(c'))\\
&=\lambda(\widetilde{f}_g(c),\widetilde{f}_h(c')).
\end{align*}
If $\f_g(c)\neq \f_h(c')$ one uses a similar argument to get the same equality. Otherwise
\[\lambda(f_g(c),f_h(c'))=\lambda(\widetilde{f}_g(c),\widetilde{f}_h(c'))=0.\]
\end{proof}

From the previous theorem one gets that same kneading sequences imply combinatorial equivalence, this is the content of the following corollary which completes the proof of item (1) of Theorem \ref{main:thm:1}.

\begin{cor}\label{cor:ce} 
Assume the same setting and hypothesis as Theorem \ref{thm:ft}. Then
$\F=\{f_1, \ldots, f_n\}$ is combinatorially equivalent to $\widetilde{\F}=\{\f_1, \ldots, \f_n\}$.
\end{cor}

\begin{proof}
Given a turning point $c$ and $g\in G_c$ define $\varphi(f_g(c))=\f_g(c)$.
From Theorem \ref{thm:ft}, it follows that $\varphi$ is well defined and 
$\lambda(f_{g_1}(c_1),f_{g_2}(c_2))=\lambda(\f_{g_1}(c_1),\f_{g_2}(c_2))$,
for all turning points $c_1,c_2$ and $g_1\in G_{c_1}$, $g_2\in G_{c_2}$. Hence,
$\varphi$ is an order preserving bijection from $\cc_G$ into $\widetilde{\cc}_G$. Now, let $x\in \textrm{Dom}(f_g)\cap\cc_G$.
So there is a turning point $c$ and $h\in G_c$ such that $x=f_h(c)$.
Note that $G_{f_h(c)}=\widetilde{G}_{\f_h(c)}$. Thus, by definition of $\varphi$ we have
\[
(\varphi\circ f_g)(x)=(\varphi\circ f_g)(f_h(c))=\varphi(f_{hg}(c))=\f_{hg}(c)=
\f_g(\f_h(c))=(\f_g\circ\varphi)(x).
\]
This concludes the proof of the corollary.
\end{proof}

The next example shows that, in general, we can not extend the topological conjugacy guaranteed by Theorem \ref{main:thm:1} into a conjugacy between invariant sets strictly bigger than the closure of the critical orbits.

\begin{example}
Let $\F=\{f_1,f_2\}$, $\widetilde{\F}=\{\widetilde{f}_1,\widetilde{f}_2\}$ be given by ${\rm Dom} (f_1)={\rm Dom} (\widetilde{f}_1)=[-1,0]$, ${\rm Dom} (f_2)={\rm Dom} (\widetilde{f}_2)=[0,1]$ and
\begin{align*}
&f_1(x)=-x^2,\, f_2(x)=-x^2,\\
&\widetilde{f}_1(x)=-\sqrt{-x}, \,  \widetilde{f}_2(x)=-\sqrt{x}.
\end{align*}
Observe that $\cc_G=\{-1,0,1\}=\widetilde{\cc}_G$, and hence $\F$ in $\cc_G$ is topologically conjugate to $\widetilde{\F}$ in $\widetilde{\cc}_G$. On the other hand, there exist no $K,\widetilde{K}$, $\F,\widetilde{\F}$-invariant, containing strictly $\cc_G,\widetilde{\cc}_G$, respectively, such that  $\F$ in $K$ is topologically conjugate to the system $\widetilde{\F}$ in $\widetilde{K}$. 
\end{example}

We will now go from combinatorial equivalence to topological conjugacy, for this we will use the hypothesis that if $(x,y)$ is a connected component of $\R\setminus {\rm cl}(\cc_G)$ then $x,\, y \in \cc_G$, and the same for $\widetilde{\cc}_G$. We notice that the argument we use to finish the prove of Theorem \ref{main:thm:1} can be applied in general to extend any strictly increasing bijection, provided the sets involved satisfy the aforementioned hypothesis.

\begin{proof}[Proof of the item (2) of Theorem \ref{main:thm:1}.]
Because of Proposition \ref{prop:top} we only need to prove the reverse implication.
Let $\varphi:\cc_G\to\widetilde{\cc}_G$ given by Corollary \ref{cor:ce}.
We will show that $\varphi$ can be extended continuously to the closure of $\cc_G$.
For that, let $x\in{\rm cl}(\cc_G)\setminus\cc_G$. The point $x$ is accumulated on both sides by points in $\cc_G$, otherwise there would be an open interval of the form $(x,a)$ or $(a,x)$ contained in $\R\setminus {\rm cl}(\cc_G)$. Thus there exist $x_m^j\in\cc_G$, $j=1,2$, such that $x_m^1\uparrow x$
and $x_m^2\downarrow x$. Since $\varphi$ is  
strictly increasing there are unique
\[
\varphi_1(x):=\displaystyle\lim_{m\to\infty} \varphi(x_m^1) \ \ \textrm{and} \ \
\varphi_2(x):=\displaystyle\lim_{m\to\infty} \varphi(x_m^2). 
\]
Note that $\varphi_1(x)\leq\varphi_2(x)$. Besides that, $\varphi_j(x)$ does not depend of the
sequence $x^j_n$ which converges to $x$. If $\varphi_1(x)< \varphi_2(x)$ then $(\varphi_1(x),\varphi_2(x))$ would be a connected component of $\R\setminus {\rm cl}(\widetilde{\cc}_G)$, this would imply that $\varphi_1(x),\, \varphi_2(x)$ are in $\widetilde{\cc}_G$. This can only happen if $x\in \cc_G$. Hence we conclude that $\varphi_1(x)=\varphi_2(x)$.
Therefore, for $x\in{\rm cl}(\cc_G)\setminus\cc_G$ we can define
\[
\varphi(x):=\displaystyle\lim_{m\to\infty} \varphi(x_m),
\]
where $x_m\in\cc_G$ is some (any) sequence such that $\lim_{m\to\infty}x_m=x$. In case $x\in \cc_G$ we have three possibilities: If $x$ is isolated from other points in $\cc_G$, then $\varphi$ is clearly continuous at $x$. If $x$ is accumulated from both sides by points in $\cc_G$ then $\varphi(x)$ is also accumulated from both sides by points in $\widetilde{\cc}_G$, otherwise there would be a connected component of $\R\setminus {\rm cl}(\widetilde{\cc}_G)$ with a border point not in $\widetilde{\cc}_G$, then $\varphi$ is continuous at $x$. If $x$ is accumulated only from the left (or only from the right) by points in $\cc_G$ then $\varphi(x)$ is also accumulated from the left (or from the right) by points in $\widetilde{\cc}_G$ and $\varphi$ is continuous at $x$. We conclude that the function $\varphi$ extends to a continuous monotonous function from $\text{cl}(\cc_G)$ to $\text{cl}(\widetilde{\cc}_G)$. Since we could repeat the same argument going from $\text{cl}(\widetilde{\cc}_G)$ to $\text{cl}(\cc_G)$ one concludes that $\varphi$ is a homeomorphism.

Let $x\in \text{Dom}(f_g)\cap \text{cl}(\cc_G)$. If $x\in \cc_G$, from the definition of $\varphi$ one has $\varphi(x)\in \text{Dom}(\widetilde{f}_g)$ and $\varphi(f_g(x))=\widetilde{f}_g(\varphi(x))$. If $x\in \text{cl}(\cc_G)\setminus \cc_G$ then $x$ is accumulated on both sides by elements on $\cc_G$, we can choose sequences in $\cc_G$ such that $x_m^1\uparrow x$ and $x_m^2\downarrow x$. Since the border points of $\text{Dom}(f_g)$ are in $\cc_G$ then we can assume $x_m^j\in \text{Dom}(f_g)$ and by continuity one gets $\varphi(x)\in \text{Dom}(\widetilde{f}_g)$ and $\varphi(f_g(x))=\widetilde{f}_g(\varphi(x))$.
Following precisely the same kind of arguments above we see that $\text{cl}(\cc_G),\text{cl}(\widetilde{\cc}_G)$ are $\F,\widetilde{\F}$-invariant, respectively. 
\end{proof}
%


\section{Kneading matrix and linearizing measure}\label{MTKT}

In this section we will extend some of the results of \cite{MT} to our context, that is systems $\F=\{f_1, \ldots, f_n\}$ of strictly monotonous functions on the real line.
In the previous section we defined the address of a point in an implicit way, we will now give an explicit definition.
For that, we take the partition of $\R$ into closed
connected sets determined by the turning points
$c_i^{\ell,r}$. Notice that this is not a partition
in strict sense, two neighboring elements intersect
in a point. Denote the elements of this partition by
$P_0,\ldots,P_{l+1}$, we use the indexing such that
$P_i$ is to the left of $P_j$ iff $i<j$, thus $P_0$
and $P_{l+1}$ are unbounded closed connected sets and
$P_1,\ldots,P_l$ are closed intervals whose border
points are turning points. We will rename the turning
points as $c_1,\ldots,c_{l+1}$ such that $c_i$ is the
left border point of $P_i$. The function $D$ from
previous section will be given by
$$
D(x)= \left\{
\begin{array}{ll}\vspace{.1cm}
P_i    & \textrm{ if }  x\in P_i \textrm{ and it is not a turning point, }    \\
\dfrac{P_{i-1}+P_i}{2}  & \textrm{ if }  x=c_i. 
\end{array}
\right.
$$

We will denote by ${\bf V}$ the vector space over $\Q$ generated by the formal vectors $P_0,\ldots,P_{l+1}$. To each $x\in \R$ we assign an element of ${\bf V}[[t]]$, i.e. a formal power series with coefficients in ${\bf V}$, this is called the {\em invariant coordinate} of $x$ and it is denoted by $\theta(x)$. It is given by
\[\theta(x)=\sum_{g\in\iS_x} \sigma(f_g)D(f_g(x))t^{|g|},\]
where $\iS_x=\{g\in S_x: \text{int}({\rm Dom}(f_g))\neq \emptyset\}$.
We endow ${\bf V}[[t]]$ with a topology given by the basic sets $v+t^m {\bf V}[[t]]$, for any $v\in {\bf V}[[t]]$ and $m\in \N$. On the other hand, we put on $\bf{V}$ a translation invariant total order such that
$P_0<P_1<\cdots<P_{l+1}$, and extend it lexicographically to $\bf{V}[[t]]$. This order implies the following 

\begin{lem}\label{lem:mon}
If $x<y$ then $\theta(x)\leq \theta(y)$.
\end{lem}
\begin{proof}
Choose $x<y$. Let $m$ be the biggest number such that $\{g\in\iS_x: |g|\leq m\}=\{g\in\iS_y: |g|\leq m\}$. If $m=\infty$ then $\iS_x=\iS_y$, and for all $g\in\iS_x$ we have $\sigma(f_g)f_g(x)<\sigma(f_g) f_g(y)$, which implies $\sigma(f_g)D(f_g(x))\leq \sigma(f_g) D(f_g(y))$ and thus $\theta(x)\leq \theta(y)$. When $m<\infty$ we still have $\sigma(f_g)D(f_g(x))\leq \sigma(f_g) D(f_g(y))$ for all $g\in\iS_x$, $|g|\leq m$. The maximality of $m$ implies that for some $g$, with $|g|=m$, we have $\sigma(f_g)D(f_g(x))< \sigma(f_g) D(f_g(y))$ and therefore $\theta(x)<\theta(y)$.
\end{proof}

The invariant coordinate $\theta(x)$ can be written as
\[\theta(x)=\theta_0(x) P_0+\cdots+\theta_{l+1}(x) P_{l+1},\]
where $\theta_j \in \Q[[t]]$, for $j=0,\ldots,l+1$. Given $P_j$, consider the generators $a_{i_1},\ldots,a_{i_k}$ such that $P_j\subset \text{Dom}(f_{a_{i_m}})$, $m=1,\ldots,k$. Define
\[e_j(t)=1-\sigma(f_{a_{i_1}})t-\cdots-\sigma(f_{a_{i_k}})t.\]
For example $e_0(t)=e_{l+1}(t)=1$.

\begin{lem}\label{lem:ld}
Let $x\in \R$ such that $f_g(x)\neq c_i$, for any $i$ and $g\in\iS_x$, then we have
\[
\sum_{j=0}^{l+1} \theta_j(x) e_j(t)=1.
\]
\end{lem}
\begin{proof}
Since ${\bf V}[[t]]$ is a free module over $\Q[[t]]$ with basis $P_0,\ldots,P_{l+1}$, we consider the $\Q[[t]]$--linear homomorphism $h:{\bf V}[[t]]\to \Q[[t]]$ defined by $h(P_j)=e_j(t)$. Suppose that, for some $g\in\iS_x$, we have $D(f_g(x))=P_j$ and $P_j\subset\text{Dom}(f_{a_{i_m}})$ for the generators $a_{i_1},\ldots,a_{i_k}$, then
\begin{align*}
h(\sigma(f_g) D(f_g(x))t^{|g|}) &=\sigma(f_g) (1-\sigma(f_{a_{i_1}})t-\cdots-\sigma(f_{a_{i_k}})t)t^{|g|}\\
&=\sigma(f_g)t^{|g|}-\sigma(f_{ga_{i_1}})t^{|g|+1}-\cdots-\sigma(f_{ga_{i_k}})t^{|g|+1}.
\end{align*}
Hence
\begin{align*}
\sum_{\substack{g\in \iS_x\\|g|=m}} h\left(\sigma(f_g) D(f_g(x))t^{|g|}\right)&= \sum_{\substack{g\in \iS_x\\|g|=m}} \left( \sigma(f_g)t^{|g|}\,\, - \sum_{\substack{a_i\\ga_i\in \iS_x}} \sigma(f_{ga_{i}})t^{|g|+1} \right)\\
&=\sum_{\substack{g\in \iS_x\\|g|=m}} \sigma(f_g)t^{|g|}\,\, -\sum_{\substack{g\in \iS_x\\|g|=m+1}} \sigma(f_g)t^{|g|}.
\end{align*}
Summing over all $m$ we get $h(\theta(x))=1$, which proves the desired equality.
\end{proof}

We say that a pair $(g,x)\in \iS_x\times \R$ is pre-turning if $f_g(x)$ is a turning point and $f_h(x)$ is not a turning point for all $h\neq g$ in the geodesic connecting $e$ with $g$. We say that a point $x\in \R$ is pre-turning if for some $g\in \iS_x$ the pair $(g,x)$ is pre-turning, in this case $x$ is said a pre-turning point of order $|g|$. 
For a fixed $m$, the set of pre-turning points of order less or equal than $m$ is finite, consider the complement of this set on $\R$, it is easy to see that if $x$ and $y$ are on the same connected component of this complement, then $\theta(x)-\theta(y)\in t^{m+1} {\bf V}[[t]]$. This implies two things: if $x_0$ is not a pre-turning point then the function $x\mapsto \theta(x)$ is continuous at $x_0$. If $x_0$ is a pre-turning point, lateral limits exist but they are not equal. We denote those limits by $\theta(x_0^-)$ (left limit) and $\theta(x_0^+)$ (right limit), that is
\[
\theta(x_0^-)=\lim_{\substack{y\to x_0\\y<x_0}}\theta(y) \ \ \textrm{and} \ \ \theta(x_0^+)=\lim_{\substack{y\to x_0\\x_0<y}}\theta(y).
\]

The {\em kneading increments} are defined by $\vt_i = \theta(c_i^+)-\theta(c_i^-)$ for $i=1,\ldots,l+1$ which measures the discontinuity of $\theta$ at $c_i$. We can write the kneading increments in terms of the basis $P_i$ as
\[\vt_i= N_{i,0}P_0+\cdots+N_{i, l+1} P_{l+1}.\]
The $(l+1)\times (l+2)$ matrix $\{N_{i,j}\}$ is called {\em kneading matrix}, see Remark \ref{remarkkm}.

Lemma \ref{lem:ld} will allow us to find a linear relation between the columns of the kneading matrix. Notice that each $x_0$ can be approach, on the left and right, by numbers which are not pre-turning points. Thus we can use the lemma to conclude that
\[\sum_{j=0}^{l+1} \theta_j(x_0^{\pm})e_j(t)=1.\]

Using the above formula at $x_0=c_i$ and subtracting the expressions we get $\sum_{j=0}^{l+1} N_{i,j} e_j(t)=0$. Denote the columns of $\{N_{i,j}\}$ by $\Gamma_0,\ldots,\Gamma_{l+1}$, then we have proven $$\sum_{j=0}^{l+1}e_j(t) \Gamma_j=0.$$
This linear relation make it possible to define a ``determinant'' for the kneading matrix.

\begin{definition}
Define $D_i={\rm det} [\Gamma_0,\ldots,\widehat{\Gamma}_i,\ldots,\Gamma_{l+1}]$, where $\widehat{\Gamma}_i$ means we deleted the column $\Gamma_i$. We define the {\em kneading determinant} as the power series $D\in \Q[[t]]$ given by
\[D=\frac{(-1)^i}{e_i(t)} D_i,\]
it does not depend on $i$.
\end{definition}

To give one very simple example consider the system $\F$ of Example
\ref{exa:nsp} which has kneading matrix
\begin{align*}
\{N_{i,j}\} = 
\begin{bmatrix*}[c]
-1+2t & \phantom{-}1 &  0    \\
  0   &      -1      & 1-2t 
\end{bmatrix*}    
\end{align*}
with kneading determinant
\[
D=1-2t.
\]

Now we will examine the difference between lateral limits for a pre-turning point. 
For that, we consider the set
\[
\iS_{x^+}=\biggl\{g\in\iS_x:  
\begin{array}{c}
\textrm{ either } x\in{\rm int}({\rm Dom}(f_g))\\
\textrm{ or } x \textrm{ is the left border point of } {\rm Dom}(f_g)
\end{array}
\biggr\}.
\]
One can show that for all $m\geq 1$ there is $\delta>0$ such that
\[
\iS_{x^+}\cap\{|g|<m\}=\iS_{y}\cap\{|g|<m\}
\]
for all $y\in(x,x+\delta)$.
We define $\iS_{x^-}$ similarly and consider the limits 
\[
D(f_g(x^+)):=\lim_{y\to x,\,x<y}D(f_g(y)) \textrm{ and }  D(f_g(x^-)):=\lim_{y\to x,\,y<x}D(f_g(y))
\]
for $g\in\iS_{x^\pm}$. Observe that if $(g,x)\in\iS_{x^\pm}\times\R$ is pre-turning then $D(f_g(x^\pm))\neq D(f_g(x))$.

\begin{lem}\label{lem:prejump}
For any $x$ pre-turning point
\[\theta(x^+)-\theta(x^-)=\sum_{c_i} \sum_{\{g\,:\, (g,x) \textrm{ is } c_i\textrm{ pre.}\}} \vt_i t^{|g|},\]
where the second sum is over all $g$'s such that $(g,x)$ is pre-turning and $f_g(x)=c_i$.
\end{lem}

\begin{proof}
Let $d\in\iS_{x^+}$ (resp. $d\in\iS_{x^-}$), if for all $g$ in the geodesic connecting $e$ with $d$ we have that $f_g(x)$ is not a turning point, then $x\in \textrm{int(Dom}(f_d))$ and $\sigma(f_d) D(f_d(x^+))=\sigma(f_d) D(f_d(x))$ (resp. $\sigma(f_d) D(f_d(x^-))=\sigma(f_d) D(f_d(x))$). Otherwise, we consider the $g$, in the geodesic connecting $e$ and $d$, closest to $e$ and such that $f_g(x)$ is a turning point. In this case the pair $(g,x)$ is pre-turning. Let $f_g(x)=c_i$ and write $d=gh$ with $|d|=|g|+|h|$. We get
\begin{align*}
\sigma(f_d)D(f_d(x^+))&=\sigma(f_g)\sigma(f_{h})D(f_{h}(c_i^{\sigma(f_g)}))\\
\text{(resp. } \sigma(f_d)D(f_d(x^-))&=\sigma(f_g)\sigma(f_{h})D(f_{h}(c_i^{-\sigma(f_g)}))\text{)}.
\end{align*}
From this we obtain
\begin{align*}
\theta(x^+)-\theta(x^-)&= \sum_{d\in\iS_{x^+}} \sigma(f_d) D(f_d(x^+)) t^{|d|}- \sum_{\wh d\in \iS_{x^-}} \sigma(f_{\wh d}) D(f_{\wh d}(x^-)) t^{|\wh d|}\\
&=\sum_{c_i} \sum_{\{g\,:\, (g,x) \textrm{ is } c_i\textrm{ pre.}\}} \sigma(f_g) t^{|g|} \left( \sum_{h\in \iS_{c_i^{\sigma(f_g)}}} \sigma(f_h) D(f_h(c_i^{\sigma(f_g)})) t^{|h|} \right.\\ 
& \hspace{4cm} \left. - \sum_{\wh h\in \iS_{c_i^{-\sigma(f_g)}}} \sigma(f_{\wh h}) D(f_{\wh h}(c_i^{-\sigma(f_g)})) t^{|\wh h|}  \right) \\
&=\sum_{c_i} \sum_{\{g\,:\, (g,x) \textrm{ is } c_i\textrm{ pre.}\}} \vt_i t^{|g|}.&
\end{align*}
\end{proof}

Given a set $J\subset \R$ and $m\geq0$ we define the pre-turning counting sequences
\[\gamma_{i,m}(J)=\# \{(g,x)\in\iS_x\times J: |g|=m,\, (g,x)\text{ pre-turning with } f_g(x)=c_i\}.
\]
Associated to this counting sequences we have the formal power series
\[\gamma_i(J)= \sum_{m\geq 0} \gamma_{i,m}(J) t^m.\]

We use Lemma \ref{lem:prejump} to prove that the ``growth'' of the function $\theta$ is determined by the pre-turning counting series and the kneading matrix.

\begin{prop}\label{lem:gteta}
Let $a<b$ and $J=[a,b]$, then
\[\theta(b^+)-\theta(a^-)=\sum_{c_i}\vt_i \gamma_i(J).\]
\end{prop}
\begin{proof}
Firstly, observe that 
$$\theta(b^+)-\theta(a^-)= \sum_{a\leq x\leq b \text{ pre-tur. }} \theta(x^+)-\theta(x^-).$$
From Lemma \ref{lem:prejump} we get
\begin{align*}
\theta(b^+)-\theta(a^-) & = \sum_{{a\leq x\leq b}\atop {\text{ pre-tur.}}}\sum_{c_i}\sum_{{(g,x) \textrm{ pre-tur. }} \atop {f_g(x)=c_i}}\vt_i t^{|g|}=\sum_{c_i}\sum_{m=0}^{\infty}\sum_{{{a\leq x\leq b} \atop {(g,x) \textrm{ pre-tur. }}} \atop {f_g(x)=c_i,\, |g|=m}}
\vt_i t^{|g|}\\
& = \sum_{c_i}\sum_{m=0}^{\infty}\vt_it^m\gamma_{i,m}(J)=\sum_{c_i}\vt_i\gamma_i(J).
\end{align*}
\end{proof}

Now we are going to investigate how $\vt_i$ relates to $\theta(c_i)$. For the next lemma we will use the set $\iS'_{c_i}=\{g\in\iS_{c_i}: c_i\in \partial {\rm Dom}(f_g)\}$, where $\partial$ refers to the boundary of the set.

\begin{lem}\label{lem:cipm}
Given a turning point $c_i$, we have
$\iS_{c_i}\setminus \iS'_{c_i} = \iS_{c_i^+}\cap \iS_{c_i^-}$. 
Moreover, $\iS_{c_i}$ decomposes as the disjoint union of the sets $\iS_{c_i^+}\cap \iS_{c_i^-}$, $\iS_{c_i^+}\cap \iS'_{c_i}$ and $\iS_{c_i^-}\cap \iS'_{c_i}$.
\end{lem}

\begin{proof}
If $g\in\iS_{c_i}\setminus \iS'_{c_i}$ then $c_i\in {\rm int}({\rm Dom}(f_g))$ and from this $g\in \iS_{c_i^+}\cap \iS_{c_i^-}$. On the other hand, if $h\in \iS_{c_i^+}\cap \iS_{c_i^-}$ then there exists $\delta>0$ such that $(c_i-\delta,c_i+\delta)\subset{\rm Dom}(f_h)$ and so $c_i\in {\rm int}({\rm Dom}(f_g))$, which implies $h\in \iS_{c_i}\setminus \iS'_{c_i}$. This proves $\iS_{c_i}\setminus \iS'_{c_i} = \iS_{c_i^+}\cap \iS_{c_i^-}$. The decomposition of $\iS_{c_i}$ follows easily from the previous equalities.
\end{proof}

Using Lemma \ref{lem:cipm} we can rewrite $\vt_i$ in the following way
\begin{align*}
\vt_i=\sum_{c_j} \sum_{\{h\in \iS_{c_i}\setminus \iS'_{c_i}: f_h(c_i)=c_j\}} (P_j-P_{j-1})t^{|h|}&+
\sum_{g\in \iS_{c_i^+}\cap \iS'_{c_i}}\sigma(f_g)D(f_g(c_i^+))t^{|g|}\\
&-\sum_{g\in \iS_{c_i^-}\cap \iS'_{c_i}}\sigma(f_g)D(f_g(c_i^-))t^{|g|}.
\end{align*}

Note that for every $h\in \iS_{c_i}\setminus \iS'_{c_i}$ such that $f_h(c_i)$ is a turning point, there is a generator $a_k$ such that $ha_k\in \iS'_{c_i}$. Thus $$\#\{h\in \iS_{c_i}\setminus \iS'_{c_i}: f_h(c_i)=c_j,\,|h|=m\}\leq \#\{g\in \iS'_{c_i}:\, |g|=m+1\}.$$
From this we define $\ell'_{c_i}(m)= \#\{g\in \iS'_{c_i}:\, |g|=m\}$ and its associated series $L'_{c_i}=\sum_{m\geq 1} \ell'_{c_i}(m)t^{m-1}$.
The important conclusion obtained from the previous discussion is that the radius of convergence of the power series $N_{i,j}$ is greater or equal to the radius of convergence of $L'_{c_i}$. 

Given an integer $m\geq1$, we are going to consider triples $(h,c_i,g)$ such that $h\in \iS^{-1}\cap G_{c_i}$, $g\in \iS'_{c_i}$, $|g|+|h|=m$ and $(h^{-1},f_h(c_i))$ is pre-turning, where $\iS^{-1}=\{t\in S^{-1}: \text{int}({\rm Dom}(f_t))\neq \emptyset\}$. Denote the set of such triples by ${\rm Tri}(m)$. Analogously we shall also consider the set $\iS=\{g\in S: \text{int}({\rm Dom}(f_g))\neq \emptyset\}$.

\begin{lem}\label{lem:bibo}
There is a bijection between ${\rm Tri}(m)$ and the set
\[\{ (d,x)\in \iS\times \R: |d|=m,\, x\in \partial {\rm Dom}(f_d)\}.\]
\end{lem}
\begin{proof}
Let $(h,c_i,g)$ be a triple, we assign to it $(h^{-1}g, f_{h}(c_i))$. It is clear that $|h^{-1}g|=m$. 
By hypothesis $c_i\in {\rm Im}(f_{h^{-1}})\cap\partial {\rm Dom}(f_g)$. This implies that either $g \notin S_{c_i^+}$ or $g\notin S_{c_i^-}$, and therefore $f_h(c_i)\in \partial {\rm Dom}(f_{h^{-1}g})$. To get ${\rm int}({\rm Dom}(f_{h^{-1}g}))\neq\emptyset$ it is enough to prove that $c_i\in{\rm int}({\rm Dom}(f_h))$. For that, consider the geodesic $g_0=e\rightarrow g_1\rightarrow\cdots\rightarrow g_{k-1}\rightarrow g_k=h^{-1}$ connecting $e$ with $h^{-1}$. We will prove that $f_{g_m}(f_h(c_i))\in{\rm int}({\rm Dom}(f_{g_m^{-1}}))$ for $m=0,1,\ldots,k$. It is immediate that $f_h(c_i)\in{\rm int}({\rm Dom}(f_e))$. Suppose that the assertion holds for $m\leq k-1$ and observe that $g_{m+1}=g_ma$ for some generating element $a$. Since $(h^{-1}, f_h(c_i))$ is pre-turning point, it follows that $f_{g_m}(f_h(c_i))\in{\rm int}({\rm Im}(f_{g_m})\cap{\rm Dom}(f_a))$ and hence $f_{g_{m+1}}(f_h(c_i))\in{\rm int}({\rm Im}(f_{g_{m+1}}))={\rm int}({\rm Dom}(f_{g^{-1}_{m+1}}))$. Therefore, $c_i=f_{h^{-1}}(f_h(c_i))\in{\rm int}({\rm Dom}(f_h))$.

Now take $x\in \partial {\rm Dom}(f_d)$ with $d\in\iS$, consider the geodesic going from $e$ to $d$ and let $b$ be the element of this geodesic which is closest to $e$ and such that $f_b(x)$ is a turning point. Such $b$ exists because otherwise $x$ would not be a border point of ${\rm Dom}(f_d)$. To the pair $(d,x)$ we assign the triple $(b^{-1},f_b(x), b^{-1}d)$. Note that $d\in\iS$ implies that $b^{-1}d\in\iS'_{f_b(x)}$ and $b\in\iS_x$.

It is not difficult to see that the assignments from $\text{Tri}(m)$ to $\{ (d,x)\in \iS\times \R: |d|=m,\, x\in \partial {\rm Dom}(f_d)\}$, and the other way around, are one inverse to the other. 
Therefore, this two sets have the same cardinality.
\end{proof}

Now we will consider the sequences $$\ell(m|J)=\#\{g\in \iS: |g|=m,\, {\rm int (Dom}(f_g)\cap J)\neq \emptyset\},$$ where $J$ is a subset of $\R$, and $\ell_x(m)=\#\{g\in \iS_x: |g|=m,\, x\in {\rm int( Dom}(f_g))\}$ for $x\in \R$. Associated to these sequences we have the formal power series
\[L(J)=\sum_{m\geq 1} \ell(m|J) t^{m-1}\, \text{ and }\, L_x=\sum_{m\geq 1} \ell_x(m) t^{m-1}.\]
When $J=\R$ we will write $\ell(m)$ and $L$ instead of $\ell(m|\R)$ and $L(\R)$.
Using Lemma \ref{lem:bibo} we get
\begin{align*}
\ell(m)&= \sum_{k=0}^{m-1} \sum_{c_i} \sum_{\{|h|=k:\, (h^{-1},f_h(c_i))\, \text{pre.}\}} \sum_{\{|g|=m-k:\, g\in \iS'_{c_i}\}} \tfrac{1}{2}\\
&=\tfrac{1}{2}\sum_{k=0}^{m-1} \sum_{c_i} \ell'_{c_i}(m-k) \sum_{\{|h|=k:\, (h^{-1},f_h(c_i))\, \text{pre.}\}} 1\\
&=\tfrac{1}{2}\sum_{k=0}^{m-1} \sum_{c_i} \ell'_{c_i}(m-k) \gamma_{i,k},
\end{align*}
and hence
\begin{align*}
L=\sum_{m\geq 1}\ell(m)t^{m-1} &= \tfrac{1}{2}\sum_{c_i}\sum_{m\geq 1} \left[ \sum_{k=0}^{m-1}\ell'_{c_i}(m-k) \gamma_{i,k} t^{m-1}\right]\\
&=\tfrac{1}{2}\sum_{c_i}\left(\sum_{m\geq 1}\ell'_{c_i}(m)t^{m-1}\right) \left(\sum_{k\geq 0}\gamma_{i,k}t^k\right)\\
&=\tfrac{1}{2}\sum_{c_i}L'_{c_i}\gamma_i.
\end{align*}

\begin{thm}\label{thm:ss0}
Let 
\[s=\limsup_{m\to \infty} \ell(m)^{1/m} \text{ and } s_0= \max_{c_i} \limsup_{m\to \infty} \ell'_{c_i}(m)^{1/m}.\]
If $s_0<s$ then $L$ extends to a meromorphic function on $\{t\in \C:|t|<1/s_0\}$ with a pole at $t=1/s$. Moreover, the kneading determinant $D$ vanishes at $t=1/s$ and it has no zero in the open disk $|t|<1/s$.
\end{thm}
\begin{proof}
By Proposition \ref{lem:gteta} we get
\[P_{l+1}-P_0=\theta(c_{l+1}^+)-\theta(c_1^-)=\sum_{c_i} \vt_i \gamma_i= \sum_{j=0}^{l+1} \left(\sum_{i=1}^{l+1}N_{i,j}\gamma_i\right)P_j.\]
This equation can be seen as a matrix multiplication, thus we can invert it, deleting one column from the kneading matrix. Since $N_{i,j}$ defines an holomorphic function for $|t|<1/s_0$, then $\gamma_i$ extends to a meromorphic function on $|t|<1/s_0$. Now, equation $L=\tfrac{1}{2}\sum_{c_i} L'_{c_i}\gamma_i$ together with the fact that $L'_{c_i}$ is holomorphic in $|t|<1/s_0$ allow us to conclude that $L$ extends to a meromorphic function on $\{t\in \C:|t|<1/s_0\}$.

Note that $1/s$ is the radius of convergence of $L$, then $L$ cannot be extended analytically in a bigger disk than $|t|<1/s$. This implies that the meromorphic extension of $L$ must have a pole in the circle $|t|=1/s$. Now, given that the coefficients of the series $L$ are nonnegative real numbers we conclude that $L$ has a pole at $t=1/s$. The existence of such pole immediately implies $D$ has a zero at $t=1/s$, otherwise $\gamma_i$, and thus $L$, would be holomorphic around $t=1/s$.

Now we will proof that in fact $D$ has no zero on the disk $\{t\in \C: |t|<1/s\}$. Choose non pre-turning points $x_j\in P_j$, $j=0,\ldots,l+1$. Assume that $x_0$ is such that $x_0<y$ for all $y\in \bigcup_{i=1}^{n} {\rm Im}(f_i)$. For each $j$, let $A_j$ be the set of generators $a_i$ such that $x_j\in {\rm Dom}(f_{a_i})$. We have
\begin{align*}
\theta(x_j)-\theta(x_0)&=P_j-P_0+\sum_{g\in A_j} \sigma(f_{g})\theta(f_{g}(x_j))t\\
&= P_j-P_0+\sum_{g\in A_j} \sigma(f_{g})\theta(x_0)t\\
& \hspace{1.7cm}+\sum_{g\in A_j} \sigma(f_{g})[\theta(f_{g}(x_j))-\theta(x_0)]t.
\end{align*}
Let $m_j=\sum_{g\in A_j} \sigma(f_{g})$, using Lemma \ref{lem:gteta} we obtain
\[\sum_{c_i} \gamma_i ([x_0,x_j])\vt_i= P_j-P_0 + m_j P_0 t + \sum_{g\in A_j} \sigma(f_{g}) \sum_{c_i}\gamma_i([x_0,f_{g}(x_j)])\vt_i t.\]
Rewriting the equality one gets
\[P_j-P_0+m_jP_0t= \sum_{c_i} \left(\gamma_i([x_0,x_j])-\sum_{g\in A_j} \sigma(f_{g})\gamma_i([x_0,f_{g}(x_j)])t\right)\vt_i.\]
Define $\Gamma_{j,i}=\gamma_i([x_0,x_j])-\sum_{g\in A_j} \sigma(f_{g})\gamma_i([x_0,f_{g}(x_j)])t$, then \[P_j-P_0+m_jP_0t=\sum_{m=0}^{l+1}\sum_{i=1}^{l+1} \Gamma_{j,i} N_{i,m} P_m,\]
which implies
\[\delta_{j,m}= \sum_{i=1}^{l+1} \Gamma_{j,i} N_{i,m},\]
for all $1\leq m, j\leq l+1$. This proves that the kneading determinant $D$ does not vanish for $|t|<1/s$.
\end{proof}
\begin{rem}
Note that if $s=s_0$ the last part of the proof of the previous theorem still works, which implies that if $s=s_0$ then $D$ has not zero in the open disk $|t|<1/s_0$.
\end{rem}
As in \cite{MT}, supposing that $s_0<s$, we use the fact that $L$ has a pole at $t=1/s$ to find a special type of measure. For $J\subset \R$ define
\[\Lambda_0(J)= \lim_{t\to1/s} \frac{L(J)(t)}{L(t)}.\]
Note that $0\leq \Lambda_0(J) \leq 1$ and $\Lambda_0([c_1,c_{l+1}])=\Lambda_0(\R)=1$. Notice that if $g\in S$, $|g|=k+1$,  then one can write in a unique way $g=a_ig'$, such that $|g'|=k$ and $a_i$ is a generator. Moreover, ${\rm int( Dom}(f_g) \cap J ))\neq \emptyset$ if and only if ${\rm int( Dom}(f_{g'}) \cap f_{a_i}(J) ))\neq \emptyset$. This implies that
\[\ell(k+1|J)=\ell(k|f_{a_1}(J))+\cdots+\ell(k|f_{a_n}(J)),\]
and from this we get
\begin{equation}\label{eq:promedida}
\Lambda_0(J)=\tfrac{1}{s} \left[\Lambda_0(f_{a_1}(J))+\cdots+\Lambda_0(f_{a_n}(J))\right].
\end{equation}
It is not true in general that $\Lambda_0$ defines a measure (in contrast to what happens in the classical context of \cite{MT}). However, we will see that from $\Lambda_0$ we can construct an actual measure and it will have the important property in equation (\ref{eq:promedida}).

Consider the set\footnote{If $\ell_x(m)=0$ then we assume $\tfrac{1}{m} \log \ell_x(m)=-\infty < \log s$.}
\[\mathfrak{A}=\{x\in \R: \limsup_{m\to\infty} \tfrac{1}{m} \log \ell_x(m) < \log s\}.\]
Define $\mathscr{A}$ to be the set of all subsets of $\R$ which can be represented as a finite union of intervals of the form $(-\infty, x)$, $(-\infty, x]$, $(x,y)$, $[x,y)$, $(x,y]$, $[x,y]$, $[x,+\infty)$, $(x, +\infty)$ for $x,\, y \in\mathfrak{A}$. It is not difficult to prove that $\mathscr{A}$ is an algebra. Let $\mathscr{B}_{\F}$ be the $\sigma$-algebra generated by $\mathscr{A}$. With these definitions we have the following rewriting of Theorem \ref{main:thm:2}. 

\begin{thm}\label{thm:measureLamb}
Assume we are in the context of the previous discussion. There exists a measure $\Lambda:\mathscr{B}_{\F} \to \R$ such that
\[\Lambda(J)=\tfrac{1}{s} \left[\Lambda(f_{a_1}(J))+\cdots+\Lambda(f_{a_n}(J))\right]\]
for all $J\in \mathscr{B}_{\F}$, and $\Lambda(J)=\Lambda_0(J)$ for all $J\in \mathscr{A}$.
\end{thm}

\begin{proof}
We will see that $\Lambda_0$ restricted to $\mathscr{A}$ is finitely additive. Let $A_1,\ldots,A_l \in \mathscr{A}$ disjoint sets, we can write
\[U:= \bigcup_{i=1}^{l} A_i=\bigcup_{i=1}^{k} J_i,\]
where the $J_i$'s are disjoint intervals whose boundary points are in $\mathfrak{A}$. We have
\[\ell (m|U)=\{g\in S: |g|=m,\, {\rm int(Dom}(f_g)\cap U)\neq \emptyset\},\]
then $\ell(m|U)\leq \sum_{i=1}^k \ell (m|J_i)$ which implies $\Lambda_0(U)\leq \sum_{i=1}^k \Lambda_0(J_i)$. On the other hand, notice that if ${\rm int (Dom}(f_g) \cap J_i) \neq \emptyset$ and ${\rm int (Dom}(f_g) \cap J_j) \neq \emptyset$, for $i\neq j$, then $\partial J_i \cap {\rm int (Dom}(f_g)) \neq \emptyset$ (and also for  $\partial J_j$). Therefore
\[\ell(m|U)\geq \sum_{i=1}^k \ell(m|J_i)-\sum_{x\in \bigcup_{i=1}^k \partial J_i} k\cdot \ell_x(m),\]
and from this we get
\[\frac{L(U)(t)}{L(t)} \geq \sum_{i=1}^k \frac{L(J_i)(t)}{L(t)}- k\cdot \sum_{x\in \bigcup_{i=1}^k \partial J_i} \frac{L_x(t)}{L(t)}.\]
Now, taking limit as $t$ goes to $1/s$, using the fact that $\bigcup_{i=1}^k \partial J_i \subset\mathfrak{A}$ and the definition of $\mathfrak{A}$ we obtain
\[\Lambda_0(U) \geq \sum_{i=1}^k \Lambda_0(J_i).\]
We had already proven the reverse inequality, then we conclude
\[\Lambda_0(U) = \sum_{i=1}^k \Lambda_0(J_i)= \sum_{j=1}^l \sum_{J_i \subset A_j} \Lambda_0(J_i)= \sum_{j=1}^l \Lambda_0 (A_j).\]
Now assume that $A_1,\,A_2,\ldots,A_l,\ldots$ are disjoint elements in $\mathscr{A}$ such that $U=\bigcup_{i=1}^{\infty} A_i\in \mathscr{A}$. Using that $\Lambda_0$ is finitely additive and the definition of $\Lambda_0$ we get
\[\sum_{i=1}^l \Lambda_0(A_i)= \Lambda_0\left( \bigcup_{i=1}^l A_i \right)\leq \Lambda_0(U) \leq \sum_{i=1}^{\infty} \Lambda_0(A_i).\]
Making $l$ go to infinity we obtain $\Lambda_0(U)=\sum_{i=1}^{\infty} \Lambda_0(A_i)$. Thus, $\Lambda_0\restriction_{\mathscr{A}}$ is a pre-measure, i.e. it is finitely additive and $\sigma$-additive for a family whose union is still in $\mathscr{A}$. Therefore, we can use the Hahn-Kolmogorov extension theorem to conclude that there is a $\sigma$-algebra $\mathscr{B}$, containing $\mathscr{A}$, and a measure $\Lambda$ in it such that $\Lambda(J)=\Lambda_0(J)$ for all $J \in \mathscr{A}$. Let $\mathscr{B}_{\F}$ be the $\sigma$-algebra generated by $\mathscr{A}$, we clearly have $\mathscr{B}_{\F} \subset \mathscr{B}$.

Now we will prove that $\Lambda$ verifies equation (\ref{eq:promedida}) for $J\in \mathscr{B}_{\F}$. First, we will see that $\mathscr{A}$ and $\mathscr{B}_{\F}$ are $\F$ forward invariant. Indeed, consider $f_{a_j}$, for $a_j$ a generator, and $U= \bigcup_{i=1}^{k} J_i \in \mathscr{A}$ such that the $J_i$'s are disjoint intervals whose border points are in $\mathfrak{A}$. Note that
\[f_{a_j}(U)=\bigcup_{i=1}^{k} f_{a_j}(J_i \cap {\rm Dom} (f_{a_j}))\]
and every $f_{a_j}(J_i \cap {\rm Dom} (f_{a_j}))$ is an interval whose border points are contained in $f_{a_j}(\partial J_i) \cup \partial {\rm Im} (f_{a_j})$. 
The assumption $s_0 < s$ implies that $\partial {\rm Im} (f_{a_j}) \subset\mathfrak{A}$. 
From this and the definition of $\mathfrak{A}$ it follows that $f_{a_j}(\mathfrak{A})\subset\mathfrak{A}$.
Therefore, we get $f_{a_j}(\partial J_i) \cup \partial {\rm Im} (f_{a_j}) \subset\mathfrak{A}$ and then $f_{a_j}(U) \in \mathscr{A}$.
Consider now the set $$\mathscr{L}=\{X \subset \R: f_{a_j}(X) \in \mathscr{B}_{\F}\},$$ note that $\mathscr{A} \subset \mathscr{L}$. Moreover, if $X_1,\ldots,X_l,\ldots$ are elements in $\mathscr{L}$ then $f_{a_j}(\bigcup_i X_i)= \bigcup_i f_{a_j}(X_i) \in \mathscr{B}_{\F}$, which shows that $\bigcup_i X_i \in \mathscr{L}$. We also have
$f_{a_j}(\R\setminus X_1)={\rm Im}(f_{a_j})\setminus f_{a_j}(X_1),$
which is an element of $\mathscr{B}_{\F}$ since ${\rm Im}(f_{a_j}) \in \mathscr{A}$. This shows that $\mathscr{L}$ is a $\sigma$-algebra, and then $\mathscr{B}_{\F} \subset \mathscr{L}$.

Let $\varepsilon>0$ and $X \in \mathscr{B}_{\F}$ such that $\Lambda(X)<\varepsilon$. By the proof of the Hahn-Kolmogorov theorem (see \cite{Tao}) there is a family of subsets $\{E_i\}_{i=1}^{\infty}$ of $\mathscr{A}$ such that $X\subset\bigcup_i E_i$ and
$ \sum_{i=1}^{\infty} \Lambda(E_i)\leq 2\varepsilon.$
Given that $E_i\in\mathscr{A}$ then
\[ \sum_{j=1}^n \Lambda(f_{a_j}(X)) \leq \sum_{i=1}^{\infty} \sum_{j=1}^n \Lambda(f_{a_j}(E_i))= s\cdot \sum_{i=1}^{\infty} \Lambda(E_i) < 2s\varepsilon.\]

Now we can prove equation (\ref{eq:promedida}). Let $J \in \mathscr{B}_{\F}$, since this $\sigma$-algebra is generated by $\mathscr{A}$ then for every $\varepsilon>0$ there is $\widetilde{J}\in\mathscr{A}$ such that $\Lambda(J\triangle \widetilde{J})<\varepsilon$. We have
\begin{align*}
\left|\Lambda(J)- \frac{1}{s}\sum_{j=1}^n \Lambda(f_{a_j}(J)) \right| &\leq |\Lambda(J)-\Lambda(\widetilde{J})|+ \frac{1}{s}\left|\sum_{j=1}^n \left[\Lambda(f_{a_j}(\widetilde{J}))-\Lambda(f_{a_j}(J))\right] \right|\\
&\leq \Lambda(J\triangle\widetilde{J})+ \frac{1}{s}\sum_{j=1}^n \Lambda(f_{a_j}(\widetilde{J} \triangle J))
<\varepsilon+ 2\varepsilon.
\end{align*}
Since $\varepsilon$ can be chosen arbitrarily small we obtain $\Lambda(J)=\tfrac{1}{s}\sum_{j=1}^n \Lambda(f_{a_j}(J))$.
\end{proof}

As consequence we obtain the following corollary.

\begin{cor}\label{cor:nonatom}
Let $\mathfrak{A}$ as in the proof of Theorem \ref{thm:measureLamb}. Then $\Lambda(\{x\})=0$ for all $x\in\mathfrak{A}$.
\end{cor}
\begin{proof}
Since $\ell(m|\{x\})=0$ for each $m\geq1$ and $\{x\}\in \mathscr{A}$ if $x\in\mathfrak{A}$, we have $\Lambda(\{x\})=\Lambda_0(\{x\})=0$.
\end{proof}


\section{The entropy}\label{entropy}

Throughout this section we will consider the system $\F=\{f_1,\ldots,f_n\}$. Following \cites{MS1980,AM} we define topological entropy of $\F$ and show that there exists a 
continuous skew product with the same topological entropy of one. Before defining this notion let us remind the reader of the classical
entropy given by Bowen \cite{bowen-1971}.

\subsection{Bowen's entropy} \label{subsec:Bowen}
Let $(M,d)$ be a compact metric space and $f\colon M\to M$ be continuous.
The original definition of Bowen does not require compactness, but for our purpose this is enough.
Given $\ve>0$ and $m\in\N$, a set $A\subset M$ is said $(f,\ve,m)$-spanning if for all $x\in M$
there is $y\in A$ for which $d(f^j(x),f^j(y))<\ve$ for all $0\leq j<m$. A set $B\subset M$ is called
$(f,\ve,m)$-separated if for any distinct $x,y\in B$ we have $d(f^j(x),f^j(y))\geq\ve$ for some $0\leq j<m$.
We denote by $r_m(f,\ve,M)$ the smallest cardinality of an $(f,\ve,m)$-spanning set of $M$, and 
denote by $s_m(f,\ve,M)$ the largest cardinality of an $(f,\ve,m)$-separated set of $M$. Consider
\[
r(f,\ve,M)=\limsup_{m\to\infty}\tfrac{1}{m}\log r_m(f, \ve, M) 
\]
and
\[
s(f,\ve,M)=\limsup_{m\to\infty}\tfrac{1}{m}\log s_m(f, \ve, M).
\]
Finally the number 
\[
\h(f)=\lim_{\ve\to0}r(f,\ve,M)=\lim_{\ve\to0}s(f,\ve,M)
\]
is called the topological entropy of $f$. For a not necessarily $f$-invariant subset 
$N\subset M$ we denote by $h_{\text{top}}(f,N)$ the topological entropy of $f$ with respect to $N$.

\subsection{Symbolic dynamics} Most of the discussion in this subsection is an adaptation of classical 
notions. Consider the alphabet $\mathcal{A}=\{1,\ldots,n,\varnothing\}$. Let 
$\Sigma^*=\{\un\omega=(\omega_m)_{m\in\N}: \omega_m\in\mathcal{A}\setminus\{\varnothing\} 
\text{ for all } m\}$, and let $\Sigma^k=\{\un\omega=(\omega_m)_{m\in\N}: 
\omega_m\in\mathcal{A}\setminus\{\varnothing\} \text{ for } 1\leq m\leq k
\text{ and } \omega_m=\varnothing \text{ for all } m>k\}$ for all $k\geq1$. 
More $\Sigma^0=\{(\varnothing,\varnothing,\ldots,\varnothing,\ldots)\}$ consists of empty sequence. 
We shall consider the symbolic space $\Sigma:=\Sigma^*\cup\bigcup_{k\geq0}\Sigma^k$. 
Endow $\Sigma$ with the topology generated by the metric $d(\un\omega,\un\upsilon)=
\sum_{m=1}^{\infty}\frac{1}{2^m}\delta_{\omega_m,\upsilon_m}$, where
$\delta_{\omega_m,\upsilon_m}=1$ if $\omega_m\neq\upsilon_m$ and $\delta_{\omega_m,\upsilon_m}=0$ otherwise.
Then $\Sigma$ is a compact metric space. 
How the proof this is direct we leave the details to the reader. 
The {\em symbolic dynamics} is the  pair $(\Sigma,\vr)$, where $\vr\colon\Sigma\to\Sigma$
is the left shift defined by $[\vr(\un\omega)]_m=\omega_{m+1}$ for all $m\geq1$.

\subsection{The associated skew-product} 
We will use the symbolic dynamics $(\Sigma,\vr)$ defined above to build a continuous skew-product that have
the same dynamical complexity of $\F$. In particular, it follows the main result of \cite{Widodo}.
For that we need of a good definition topological entropy 
for $\F$. A classical result due to Misiurewicz and Szlenk \cite{MS1980} says that the topological
entropy of a continuous piecewise strictly monotone map is given by exponential growth rate of the 
lap numbers. More specifically, if $f\colon I\to I$ is a piecewise monotone continuous map
of a compact interval $I$ and $\ell(f^m)$ denotes the number of maximal intervals on
which $f^m$ is monotone then $\h(f)=\log s$, where $s=\lim_{n\to\infty}\ell(f^m)^{1/m}$.
Inspired by this formula and Theorem \ref{thm:ss0} we define the {\em topological entropy of $\F$ restricted to $J\subset\R$} by
\[
\h(\F;J)=\limsup_{m\to\infty}\tfrac{1}{m}\log\ell(m|J).
\]
When $J=\R$ we just write $\h(\F)$.
Even though $\F$ is far from what we call a continuous map the above formula is satisfactory,
see \cite{MiZi,Widodo,AM} where analogous formulas are used in similar contexts.
For convenience we shall use the following convention $f_\varnothing=\textrm{Id}_\R$.

Firstly,  it is easy to see that the entropy is an invariant of
topological conjugacy:

\begin{prop}\label{prop:tcentropy}
Let $\F=\{f_1, \ldots, f_n\}$ and $\widetilde{\F}=\{\f_1, \ldots, \f_n\}$ be two systems 
with ${\rm Dom}(f_i)={\rm Dom}(\f_i)$, for all $1\leq i\leq n$ and let $K,\widetilde{K}$ two invariant sets by $\F,\widetilde{\F}$, respectively.
If $\F$ in $K$ is topologically conjugate to $\widetilde{\F}$ in $\widetilde{K}$, 
then $h_{\rm top}(\F;K)=h_{\rm top}(\widetilde{\F};\widetilde{K})$.
\end{prop}
\begin{proof}
Note that $\ell(m|K)=\widetilde{\ell}(m|\widetilde{K})$ for all $m\geq1$.
\end{proof}

\noindent
{\sc Admissibility:} A sequence $\un\omega\in\Sigma$ is said {\it admissible} for $x\in\R$ 
if ${\rm int}(\textrm{Dom}(f_{\omega_m}\circ\cdots\circ f_{\omega_1}))\neq\emptyset$
and $x\in\textrm{Dom}(f_{\omega_m}\circ\cdots\circ f_{\omega_1})$ for all $m\geq1$.

\medskip

Letting $\wh I=\bigcup_{i=1}^{n}\text{Dom}(f_i)\cup\text{Im}(f_i)$ we define $\Gamma=\{(\un\omega,x)\in\Sigma\times\wh I:
\un\omega \text{ is admissible}$ 
$\text{for } x\}$.

\begin{lem}
$\Gamma$ is a compact metric space.
\end{lem}
\begin{proof}
The natural distance in $\Gamma$ is $d_\Gamma((\un\omega,x),(\un\upsilon,y))=\max\{d(\un\omega,\un\upsilon),|x-y|\}$.
Since $\Sigma\times\wh I$ is compact, to conclude that $\Gamma$ is compact it is enough to show that $\Gamma$ is closed. 
For that, let a pair $(\un\omega,x)$, and let a sequence $\{(\un\omega^l,x^l)\}_{l\in\N}$ 
in $\Gamma$ such that $\lim_{l\to\infty}(\un\omega^l,x^l)=(\un\omega,x)$.
Fixed any $m\geq1$, there exists $L\geq1$ such that 
$(\omega^l_1,\ldots,\omega^l_m)=(\omega_1,\ldots,\omega_m)$ for all $l\geq L$. Hence
${\rm int}(\textrm{Dom}(f_{\omega_m}\circ\cdots\circ f_{\omega_1}))\neq\emptyset$ and
$x^l\in\textrm{Dom}(f_{\omega_m}\circ\cdots\circ f_{\omega_1})$ for all $l\geq L$.
Therefore $x\in\textrm{Dom}(f_{\omega_m}\circ\cdots\circ f_{\omega_1})$, and so 
$(\un\omega,x)\in\Gamma$.  
\end{proof}

We consider the skew-product $F\colon\Gamma\to\Gamma$ given by
$F(\un\omega,x)=(\vr(\un\omega),f_{\omega_1}(x))$. It is easy  
to show that $F$ is continuous. Below, we state the main result of this section.

\begin{thm}\label{thm:skewentropy}
Let $\F=\{f_1,\ldots,f_n\}$ be a system and let $F$ be the skew-product associated to $\F$
defined on $\Gamma$. Then
\[
h_{{\rm top}}(F)=h_{{\rm top}}(\F).
\]
\end{thm} 

To see why this kind of result is interesting the reader could read the remark after Theorem 5.5 of \cite{Widodo}.
The proof of Theorem \ref{thm:skewentropy} will use some auxiliary lemmas. Let 
$\pi_1\colon\Sigma\times\wh I\to\Sigma$ be the projection on the first coordinate 
$\pi_1(\un\omega,x)=\un\omega$, and consider the compact $\vr$-invariant set 
$\Sigma_\Gamma:=\pi_1(\Gamma)$. We denote by $\vr_\Gamma$ the restriction
of the shift map $\vr\colon\Sigma\to\Sigma$ to $\Sigma_\Gamma$. The lemma below can
be seen as a generalization of \cite{Widodo}*{Theorem 5.5}. 

\begin{lem}\label{lem:lementropy1}
$h_{{\rm top}}(\vr_\Gamma)=h_{{\rm top}}(\F)$.
\end{lem}
\begin{proof}
Given $\un\omega\in\Sigma_\Gamma$ and $m\geq1$ we consider the cylinder
$C_m(\un\omega)=\{\un\upsilon\in\Sigma_\Gamma:\omega_i=\upsilon_i \text{ for } i=1,\ldots,m\}$.
The family of all cylinders of length $m$ is denoted by $\mathfrak{C}_m$. 
Using the definition of spanning and separated sets we have
\[
r_m(\vr_\Gamma,\tfrac{1}{2^j},\Sigma_\Gamma)\leq\#\mathfrak{C}_{m+j}
\leq s_{m+j}(\vr_\Gamma,\tfrac{1}{4},\Sigma_\Gamma)
\text{ for each } j\geq1 \text{ fixed}.
\] 
From this we get easily 
\[
h_{{\rm top}}(\vr_\Gamma)=\limsup_{m\to\infty}\tfrac{1}{m}\log\#\mathfrak{C}_m. 
\]
Since $\#\mathfrak{C}_m=\ell(1)+\cdots+\ell(m)+1$ for all $m\geq1$, the lemma follows immediately.
\end{proof}

Our next aim is to show that the entropy of the skew-product is zero on fibers.
Let $\un\omega\in\Sigma_\Gamma$. The set of points $x\in\wh I$ for which $\un\omega$ is admissible
will be denoted by $I(\un\omega)$. It is clear that $I(\un\omega)$ is compact. 

\begin{lem}\label{lem:lementropy2}
For each $\un\omega\in\Sigma_\Gamma$ holds $h_{{\rm top}}(F,\pi_1^{-1}(\un\omega))=0$.
\end{lem}
\begin{proof}
Fix $\un\omega\in\Sigma_\Gamma$.
Given $\ve>0$, we claim that there is a constant $\Upsilon_\ve$ such that 
$r_m(F,2\ve,\pi_1^{-1}(\un\omega))\leq\Upsilon_\ve m$ for every $m\geq1$.
Indeed, consider the sequence defined inductively as follows: $I_1(\un\omega):=I(\un\omega),
I_2(\un\omega):=f_{\omega_1}(I_1(\un\omega)),\ldots, 
I_l(\un\omega):=f_{\omega_{l-1}}(I_{l-1}(\un\omega)),\ldots$.
Let $\text{Co}(\un\omega)$ be the convex hull of $\bigcup_{l=1}^\infty I_l(\un\omega)$.
For each $l$ we choose a finite set $S_l\subset I_l(\un\omega)$ such that each point
in $I_l(\un\omega)$ is at most $\ve$ far from $S_l$ and $\# S_l\leq\Upsilon_\ve$, 
where $\Upsilon_\ve=\tfrac{|\text{Co}(\un\omega)|}{\ve}+1$. Fixed $m\geq1$ set
\[
S(m)=S_1\cup\bigcup_{l=1}^{m-1}\left[(f_{w_l}\circ\cdots\circ f_{w_1})^{-1}(S_{l+1})\cap I(\un\omega)\right].
\] 
Note that $\# S(m)\leq\Upsilon_\ve m$ since $f_{w_l}\circ\cdots\circ f_{w_1}\colon I_1(\un\omega)\to I_{l+1}(\un\omega)$
is a homeomorphism for all $l\geq1$. It is simple to prove that $\{\un\omega\}\times S(m)$ is a
$(F,2\ve,m)$-spanning set of $\pi_1^{-1}(\un\omega)$. Hence
the claim follows, and therefore $h_{\text{top}}(F,\pi_1^{-1}(\un\omega))=0$. 
\end{proof}

Now we can prove the main result of this section.

\begin{proof}[Proof of Theorem \ref{thm:skewentropy}.]
From Bowen's inequality \cite{bowen-1971}*{Theorem 17} we have
\[
h_{\text{top}}(\vr_\Gamma)\leq h_{\text{top}}(F)\leq h_{\text{top}}(\vr_\Gamma)+
\sup_{\un\omega\in\Sigma_\Gamma}h_{\text{top}}(F,\pi_1^{-1}(\un\omega)).
\]
The theorem is now  is a direct consequence of Lemmas \ref{lem:lementropy1} and \ref{lem:lementropy2}.
\end{proof}

\section{Examples and applications}\label{appexamples}

In this section we provide some examples and applications of the theory developed on previous sections.

\subsection{Separability Hypothesis}
Let $\F=\{f_1, \ldots, f_m, f_{m+1}, \ldots, f_{m+k}\}$ be a system such that there is a closed interval $I$ and $\beta>1$ with the following properties:
\begin{enumerate}[{\rm (i)}]
\item  $I=\bigcup_{i=1}^{m} {\rm Dom}(f_i)$ and $f_i({\rm Dom}(f_i))=I$ for all $1\leq i\leq m$.
\item $I=\bigcup_{i=m+1}^{m+k} f_i({\rm Dom}(f_i))$ and ${\rm Dom}(f_i)=I$ for all $m+1\leq i\leq m+k$.
\item $|f_i(x)-f_i(y)|\geq \beta |x-y|$ for all $1\leq i\leq m$, $x,\, y \in {\rm Dom}(f_i)$.
\item $|f_i(x)-f_i(y)|\leq \beta^{-1} |x-y|$ for all $m+1\leq i\leq m+k$, $x,\, y \in {\rm Dom}(f_i)$.
\end{enumerate}

In this case one can easily prove that $\F$ verifies the separability hypothesis and thus one can apply Theorem \ref{thm:ft} to such systems. Indeed, notice that $I=\bigcup_{i=1}^{m}f_i^{-1}(I)$ and then
\[I=\bigcup_{i_1,\ldots,i_j=1}^{m}f_{i_1}^{-1}\circ \cdots \circ f_{i_j}^{-1}(I).\]
Since the diameter of the sets $f_{i_1}^{-1}\circ\cdots\circ f_{i_j}^{-1}(I)$, for $i_1,\ldots,i_j \in [1,m]$, goes to zero as $j$ goes to infinity we conclude that ${\rm cl}(\cc_{S^{-1}})=I$. In the other hand, we have that
\[I=\bigcup_{i=m+1}^{m+k}f_i(I)=\bigcup_{i_1,\ldots,i_j=m+1}^{m+k}f_{i_j}\circ\cdots\circ f_{i_1}(I).\]
Since the diameter of the sets $f_{i_1}\circ\cdots\circ f_{i_j}(I)$, for $i_1,\ldots,i_j \in [m+1,m+k]$, goes to zero as $j$ goes to infinity we conclude that ${\rm cl}(\cc_{S})=I$. This shows that $\F$ verifies the separability hypothesis.

Now consider a system $\F=\{f_1,\ldots,f_m\}$ such that hypothesis (i) and (iii), from the previous example, are satisfied. Again we will have that ${\rm cl}(\cc_{S^{-1}})=I$, but we do not know if ${\rm cl}(\cc_{S})=I$. To every $p\in I$ we can associate a set $\Sigma_p \subset \{1,\ldots,m\}^{\N}$ defined by
\[\Sigma_p=\{(x_1,x_2,\ldots)\in \{1,\ldots,m\}^{\N} : p \in {\rm Dom}(f_{x_j}\circ\cdots\circ f_{x_1}),\, \forall j \geq 0\},\]
it is easy to prove that in fact $(x_1,x_2,\ldots,x_j,\ldots)\in \Sigma_p$ if and only if
\[\{p\} = \bigcap_{j \geq 0} f_{x_1}^{-1}\circ f_{x_2}^{-1}\circ \cdots \circ f_{x_j}^{-1}(I).\]
Observe that if $\Sigma_p$, $\Sigma_q$ contain elements coinciding in their first $j$ symbols, say $(x_1,\ldots,x_j)$, then $p,q \in f_{x_1}^{-1}\circ f_{x_2}^{-1}\circ \cdots \circ f_{x_j}^{-1}(I)$ and $|p-q|\leq \beta^{-j} |I|$. This implies that if $\Sigma_p$ contains an element of $\{1,\ldots,m\}^{\N}$ which has dense orbit with respect to the shift, then the future orbit of $p$, i.e. $\{f_g(p):\, g\in S_p\}$, is dense in $I$. We conclude that if for some turning point $c$ the set $\Sigma_c$ contains a sequence with dense orbit, then ${\rm cl}(\cc_S)=I$ and the system $\F$ would verify the separability hypothesis.

If $\F$ does not have a turning point $c$ such that $\Sigma_c$ has a sequence with dense orbit, we can modify the system to get a new one satisfying the separability hypothesis.
Firstly, from the discussion in the previous paragraph, it is not difficult to prove that the set of points $p\in I$ such that $\Sigma_p$ has a sequence with dense orbit is a dense set. Thus, taking any $1\leq i\leq m$, we can choose a point $p_0 \in {\rm Dom}(f_i)$ such that $\Sigma_{p_0}$ contains a point with dense orbit, and hence we can define the new system $\F'=\{f_1,\ldots,f_{i-1}, f_i^-,f_i^+,f_{i+1},\ldots,f_m\}$ where 
\[f_i^-=f_i\restriction_{\{q\in {\rm Dom}(f_i):\,q\leq p_0\}} \textrm{ and }\, f_i^+=f_i\restriction_{\{q\in {\rm Dom}(f_i):\,q\geq p_0\}}.\] 
From the previous discussion one sees that the new system $\F'$ verifies the separability hypothesis.

\subsection{The classical kneading sequences revisited}

Before treating the classical case we shall prove an auxiliary lemma.

\begin{lem}\label{lem:den}
Let $\F=\{f_1, \ldots, f_n\}$ and $x,y\in\R$ two different points. These points can be separated by an element in $S$ if and only if $[x,y]\cap \cc_{S^{-1}}\neq\emptyset$. 
\end{lem}
\begin{proof} 
In order to prove the first part, let $x,y\in\R$ be two different points. Assume that they can be separated by an element in $S$ and let $J\in S$ such that $x,y\in\text{Dom}(f_J)$ and $D(f_J(x))\neq D(f_J(y))$. Thus there is a turning point $c\in[f_J(x),f_J(y)]$ and hence $f_{J^{-1}}(c)\in[x,y]$. To prove the reverse implication let a turning point $c$ and $J\in S$ such that $f_{J^{-1}}(c)\in[x,y]$. 
Let $f_{s_0}, \ldots, f_{s_j}$ be a sequence so that $s_0$ is the identity of $S$, $s_j=J$ and $s_{m}=s_{m-1} a_{i(m)}$ for any $1\leq m\leq j$ and some generating element $a_{i(m)}$. We have three cases to consider: if $x,y\in\text{Dom}(f_J)$ then we are done. If there is $s_m$ such that either $x\in\text{Dom}(f_{s_m})$ and $y\notin\text{Dom}(f_{s_m})$ or $x\notin\text{Dom}(f_{s_m})$ and $y\in\text{Dom}(f_{s_m})$ then again we are done.
Lastly, let the first $s_m$ such that $x,y\in\textrm{Dom}(f_{s_{m-1}})$ and $x,y\notin\textrm{Dom}(f_{s_{m}})$. This implies that either $f_{s_{m-1}}(x)< c_{i(m)}^-<c_{i(m)}^+< f_{s_{m-1}}(y)$ or $f_{s_{m-1}}(y)< c_{i(m)}^-<c_{i(m)}^+< f_{s_{m-1}}(x)$, since $f_{J^{-1}}(c)\in\textrm{Dom}(f_{s_m})$. In any case $D(f_{s_{m-1}}(x))\neq D(f_{s_{m-1}}(y))$. 
\end{proof}

Now, let $I=[c_1, c_{l+1}]$ be a compact interval, and let $f:I\to I$ be a {\em multimodal map}. This means that $f$ is continuous and there are points, called {\em turning points}, $c_2<\cdots<c_l$ in $(c_1,c_{l+1})$ so that $f$ is strictly monotone in each interval $I_i=[c_i,c_{i+1}]$ for $i=1,\ldots,l$ and
changes monotonicity at each point $c_i$ for $i=2,\ldots,l$. Consider the alphabet $\mathcal{A}=\{c_1,I_1,c_2,\ldots,c_l,I_l,c_{l+1}\}$. For each $x\in I$ one associate a sequence $i_f(x)=(i_{f,m}(x))_{m\geq0}\in\mathcal{A}^{\{0,1,2,\ldots\}}$, called {\em itinerary}, defined as follows:
$$
i_{f,m}(x)= \left\{
\begin{array}{ll}
I_i & \textrm{ if }  x\in I_i\setminus\{c_i,c_{i+1}\}, \\
c_i  & \textrm{ if }  x=c_i. 
\end{array}
\right.
$$
The {\em kneading sequence} of $f$ is the $(l+1)$--tuple $\mathcal{K}(f)=(i_f(c_1),\ldots,i_f(c_{l+1}))$. The theorem below is well known and we will see how it follows from Theorem \ref{main:thm:1}; see \cite{dMvS} for more details on this.

\begin{thm}\label{thm:classical}
Let $f:I\to I$ and $\f:\widetilde{I}\to\widetilde{I}$ be two multimodal maps with the same number of turning points, and assume that neither of them has wandering intervals, interval of periodic points or attracting periodic points. Then $f,\f$ are topologically conjugate if and only if they have the same kneading sequences.
\end{thm}

\begin{proof}
First we observe that we can assume $c_i=\widetilde{c}_i$ for $i=1,\ldots,l+1$. Indeed, let an order preserving homeomorphism $\psi:I\to\widetilde{I}$ with $\psi(c_i)=\widetilde{c}_i$, and define $\widehat{f}:=\psi^{-1}\circ\widetilde{f}\circ\psi$. Consequently $\widehat{f}$ is a multimodal map, and $\widehat{f}$ is strictly increasing (decreasing) on $I_i$ iff $\widetilde{f}$ is strictly increasing (decreasing) on $\widetilde{I}_i$ for $i=1,\ldots,l$. Moreover, $f,\f$ are topologically conjugate iff $f,\widehat{f}$ are topologically conjugate.    

Now consider the system $\F=\{f_1,\ldots,f_l\}$ defined by $f_i=f\rr_{I_i}$ for $i=1,\ldots,l$. Observe that $\F$ separates points in the past because the intervals intersect only at border points. Letting $\cc(f)=\{x\in I :\; \exists\: m\geq0 \textrm{ such that } f^m(x)=c_i \textrm{ for some } c_i\}$ we get $\cc(f)=\cc_{S^{-1}}$. 

\medskip
\noindent
{\sc Claim:} $\textrm{cl}(\cc(f))=I$.

\medskip
The proof of this claim is standard, but we include it for completeness. Suppose that there exists an open interval $J\subset I\setminus\textrm{cl}(\cc(f))$. Since $f$ has no periodic attractors nor wandering intervals then there are $0\leq m<k$ such that $f^m(J)\cap f^k(J)\neq\emptyset$. Because $J\cap\cc(f)=\emptyset$, we may assume, without loss of
generality, that $m=0$ and $k$ is the smallest integer with such
property, and put $L=\bigcup_{j\geq0}f^{jk}(J)$.
It follows that $L$ is a non-empty interval that contains no $\{c_2,\ldots,c_l\}$ and so $f^k:L\to L$ is a strictly monotone continuous map. Thus either $L$ contains an interval of periodic points for $f^k$, or some open interval on $L$ converges to a single periodic point, which is a contradiction with the hypotheses on $f$. The proof of the claim is finished. 

\medskip
From the claim and Lemma \ref{lem:den} we have that $\F$ separates any two different points by an element in $S$, and hence $\F$ verifies the separability hypothesis. In the same way we can consider the system $\widetilde{\F}=\{\f_1,\ldots,\f_l\}$ defined by $\f_i=\f\rr_{I_i}$ and show that $\widetilde{\F}$ also satisfies the separability hypothesis. 

It is easy to see that $f,\f$ are topologically conjugate if and only if $\F$ in $I$ is topologically conjugate to $\widetilde{\F}$ in $I$. Furthermore, $f,\f$ have the same kneading sequence if and only if $\F,\widetilde{\F}$ have the same kneading sequence since $f_g(c)=f^m(c)$ for all turning point $c$ and $g\in S_c$ with $|g|=m$ (the same equality holds replacing $f$ by $\f$). Therefore, Theorem \ref{main:thm:1} implies Theorem \ref{thm:classical}.
\end{proof}

\subsection{Topological semiconjugacy to an affine map}
In this subsection we will show how the results of Section \ref{MTKT} generalize \cite{MT}.
In fact, we consider more general maps, namely piecewise continuous piecewise strictly monotone maps,
and prove that if the map has positive entropy then it is semiconjugate to a map with constant slope.
This case was already treated in \cite{Parry,AM} where their proofs do not use kneading theory. Finally, we will see how a partial version of the main result of \cite{BHV1} can be obtained using our techniques.

Let $[c_1,c_{l+1}]$ be a compact interval, and let a piecewise continuous 
piecewise strictly monotone map $f\colon[c_1,c_{l+1}]\to[c_1,c_{l+1}]$. Let $c_2<\cdots<c_l\in(c_1,c_{l+1})$ be points such that the restricions $f\restriction_{[c_i,c_{i+1}]}$ are continuous and strictly monotonous. By monotonicity, the one-sided limits exist in each
point $c_2,\ldots,c_l$, and so if $f$ is discontinuous at $c_i$ for $i=2,\ldots,l$ then we will allow that $f$ assumes both values at this point.
As before, we can also consider the associate system $\F=\{f_1,\ldots,f_l\}$ defined by 
$f_i=f\rr_{I_i}$, where $I_i=[c_i,c_{i+1}]$ for $i=1,\ldots,l$. A simple calculation shows that
the kneading matrix of $\F$ is given by

\begin{equation}\label{eq:matrix}
\{N_{i,j}\}=\left[
\begin{array}{ccc}
\begin{array}{r}
  \hspace{-.3cm}     -1 \hspace{-.4cm}\\
\end{array} & 
\begin{array}{ccc}
      N_{1,1} \hspace{.05cm} & \cdots & \hspace{.01cm} N_{1,l} \\
\end{array} & 
\begin{array}{c}
     \hspace{-.4cm}  0 \\
\end{array} \\
\begin{array}{r}
   \hspace{-.3cm}    \phantom{-}0 \hspace{-.4cm} \\ \hspace{-.3cm}\phantom{-}\vdots \hspace{-.4cm} 
	\\ \hspace{-.3cm}\phantom{-}0 \hspace{-.4cm} \\
\end{array} &     
\begin{array}{:ccc:}
        \hdashline
        N_{2,1}&\cdots& N_{2,l}\\
        \vdots&\ddots&\vdots\\
        N_{l,1}&\cdots& N_{l,l}\\
				\hdashline
      \end{array} & 
\begin{array}{c}
     \hspace{-.4cm}  0 \\ \hspace{-.4cm} \vdots \\ \hspace{-.4cm} 0 \\
\end{array} \\
\begin{array}{c}
     \hspace{-.3cm}  \phantom{-}0 \hspace{-.4cm} \\ 
\end{array} & 
\begin{array}{ccc}
       N_{l+1,1} \hspace{-.1cm} & \cdots & \hspace{-.1cm} N_{l+1,l} \\
\end{array} & 
\begin{array}{c}
      \hspace{-.4cm} 1 \\
\end{array}\\
\end{array}
\right]. 
\end{equation}

\begin{rem}\label{remarkkm}
When $f$ is a continuous piecewise strictly monotone map, then
the kneading matrix of the associate system $\F$ has the same form above
and the central sub-matrix (inside the box) is the kneading matrix of $f$ as defined by Milnor-Thurston \cite{MT}*{Section 4}.
It is not difficult to prove that $D=D_{{\rm MT}}$, where $D_{{\rm MT}}$
is the kneading determinant defined by Milnor-Thurston. This equality justifies why
\cite{MT} does not use the boundary points of the interval $[c_1, c_{l+1}]$ in its analysis.
\end{rem} 

For each $m\geq1$ and $i=1,\ldots,l+1$ we get $\ell(m)=\ell(f^m)$ and $1\leq\ell'_{c_i}(m)\leq2$, where
$\ell(f^m)$ is the lap number of $f^m$. Let $s,s_0$ is as in Theorem \ref{thm:ss0} and assume $1=s_0<s$ (this means that the entropy $\log s$ of $\F$ is positive).
Let $\mathfrak{A}$ the set defined in the proof of Theorem \ref{thm:measureLamb} and observe that $\mathfrak{A}=\R$. Thus, Theorem \ref{thm:measureLamb} and Corollary \ref{cor:nonatom} imply that there exists a non-atomic measure $\Lambda:\mathscr{B}_{\F}\to\R$, where $\mathscr{B}_{\F}$ is now the Borel $\sigma$-algebra  in $\R$.

Now, let us proceed with the construction of the semiconjugacy between $f$ and the linear model. 
Let $\varphi:[c_1,c_{l+1}]\to[0,1]$ defined by $\varphi(x)=\Lambda([c_1,x])$.
By definition it follows that $\varphi$ is a continuous, surjective and monotone map.
Note that there exist intervals (possibly degenerate) $\widetilde I_i=[\widetilde c_i,\widetilde c_{i+1}]$
with $i=1,\ldots,l$ so that $0=\widetilde c_1\leq\widetilde c_2\leq\cdots\leq\widetilde c_l
\leq\widetilde c_{l+1}=1$ and $\varphi(I_i)=\widetilde I_i$. For $y=\varphi(x)\in \widetilde I_i$, with $x\in I_i$, define $\f_i(y)=\varphi(f_i(x))$. If we also have $y=\varphi(x')$, with $x'\in I_i$, then
\[0= s \Lambda([x,x'])=\Lambda(f_i([x,x'])),\]
which implies $\varphi(f_i(x))=\varphi(f_i(x'))$. This shows that $\f_i$ is well defined. Moreover, for $y=\varphi(x)\in \widetilde{I}_i$, with $x\in I_i$, one has
\[s(\varphi(x)-\varphi(c_i))=s \Lambda([c_i,x])= \Lambda(f_i([c_i,x]))=\sigma(f_i) (\varphi(f_i(x))-\varphi(f_i(c_i))),\]
and then
\[\widetilde{f}_i(y)= \sigma(f_i)s(y-\widetilde{c}_i)+ \widetilde{f}_i(\widetilde{c}_i). \]
So the maps $\f_i$ are affine with slope $\pm s$ and the system (which may be degenerate) $\F=\{\f_1,\ldots,\f_l\}$ in $[0,1]$ is semiconjugate to $\F$ in $[c_1,c_{l+1}]$, since $\f_i\circ\varphi=\varphi\circ f_i$.

Now we will apply the previous discussion to a particular case and in this way we will obtain a partial version of the main result in \cite{BHV1}. Consider a system $\F=\{f_0,\, f_1\}$ such that ${\rm Dom}(f_0)=[0,q]$, ${\rm Dom}(f_1)=[q,1]$, for some $0<q<1$, and $0\leq f_0(0)<f_0(q)\leq 1$, $0\leq f_1(q)<f_1(1)\leq 1$. In \cite{BHV1}, the authors call the function associated to such $\F$, i.e.
\[T(x)= \begin{cases} f_0(x) & {\rm if}\,\, 0\leq x\leq q,\\
f_1(x) & {\rm if}\,\, q< x\leq 1,
\end{cases}\]
a {\em special overlapping dynamical system}.
We know that the entropy of $T$ is given by the exponential growth rate of the lap numbers, that is, $\log s$.

The system $\F$ has the turning points $c_1=0$, $c_2=q$ and $c_3=1$, they partition the real line in the intervals $P_0=(-\infty, 0]$, $P_1=[0, q]$, $P_2=[q, 1]$ and $P_3=[1,+\infty)$. Using the notation in \cite{BHV1}, we denote by $\alpha:=(\alpha_0, \alpha_1,\ldots)$ and $\beta:=(\beta_0, \beta_1,\ldots)$ two sequences of symbols (with entries in $\{0,1\}$) given by the following properties:
\begin{enumerate}[(i)]
\item $\alpha_0=0$ and $\beta_0=1$.
\item $f_{\alpha_{i-1}}\circ\cdots\circ f_{\alpha_{0}}(q) \in P_{\alpha_i +1}$ and if $f_{\alpha_{i-1}}\circ\cdots\circ f_{\alpha_{0}}(q)=q$ then $\alpha_i=0$.
\item $f_{\beta_{i-1}}\circ\cdots\circ f_{\beta_{0}}(q) \in P_{\beta_i +1}$ and if $f_{\beta_{i-1}}\circ\cdots\circ f_{\beta_{0}}(q)=q$ then $\beta_i=1$.
\end{enumerate}
It is not difficult to prove that the sequences $\alpha,\, \beta$ are related to our concepts by the equations
\[\theta(q^-)= \sum_{i\geq 0} P_{\alpha_i +1} t^i,\,\, \theta(q^+)= \sum_{i\geq 0} P_{\beta_i +1} t^i.\]
Therefore
\[\nu_2= \theta(q^+)-\theta(q^-)=\sum_{i\geq 0} (P_{\beta_i +1}-P_{\alpha_i +1}) t^i.\]
Note that $P_{x+1}=x P_2 + (1-x) P_1$ for $x \in \{0,1\}$, then we can write
\[\nu_2= \left[\sum_{i\geq 0} (\alpha_i -\beta_i) t^i \right] P_1+\left[\sum_{i\geq 0} (\beta_i -\alpha_i) t^i \right] P_2,\]
and get $N_{2,1} = \sum_{i\geq 0} (\alpha_i-\beta_i) t^i$. The main result in \cite{BHV1} is the following

\begin{stm}\label{BHV1}
Assume that $f_0(0)=0$, $f_1(1)=1$ and $f_0,\, f_1$ are expanding. Let $([0,1],T)$ be any such special overlapping dynamical system with point of discontinuity $q$, critical itineraries $\alpha,\, \beta$ and let $r$ be given by
\[r=\min \left\{t\in (0,1):\, \sum_{i=0}^{\infty} \alpha_i t^i=\sum_{i=0}^{\infty} \beta_i t^i\right\}.\]
Then:
\begin{itemize}
\item[(1)] The dynamical system $([0,1],T)$ is topologically conjugate to the piecewise affine dynamical system $([0,1],U)$, given by
\[U(x)=\begin{cases} r^{-1}x & \text{ if } 0\leq x\leq p,\\ r^{-1}(x-1)+1 & \text{ if } p< x\leq 1,\end{cases}\]
where $p=(1-r) \sum_{i=0}^{\infty} \alpha_i r^{i}$.
\item[(2)] The entropy of $([0,1],T)$ is $-\log r$.
\end{itemize}
\end{stm}

The authors of \cite{BHV1} realized that there was a problem in their argument and published an erratum in the arXiv \cite{BHV2Errata} stating that there was an error. Here we will show how to prove the following partial version of this result:

\begin{thm}
Let $\F=\{f_0,\, f_1\}$ be a system such that ${\rm Dom}(f_0)=[0,q]$, ${\rm Dom}(f_1)=[q,1]$, for some $0<q<1$, and $0\leq f_0(0)<f_0(q)\leq 1$, $0\leq f_1(q)<f_1(1)\leq 1$. Let $\alpha,\, \beta$ be the critical itineraries defined as before. If the topological entropy $\log s$ of $\F$ is positive, then
\[s^{-1}=\min \left\{t\in (0,1):\, \sum_{i=0}^{\infty} \alpha_i t^i=\sum_{i=0}^{\infty} \beta_i t^i\right\}.\]
Moreover, $\F$ in $[0,1]$ is semiconjugate to a system $U=\{U_0,U_1\}$ in $[0,1]$, where $U_0,\, U_1$ are affine maps with slope $s$. If in addition to $\log s>0$ we also assume that $f_0(0)=0$ and $f_1(1)=1$, then $U$ is of the form
\begin{align*}
U_0(x)&=sx,\,  \text{ for } 0\leq x\leq p,\\
U_1(x)&=s(x-1)+1,\,  \text{ for } p\leq x\leq 1,
\end{align*}
where $p=(1-s^{-1}) \sum_{i=0}^{\infty} \alpha_i s^{-i}$.
\end{thm}
\begin{proof}
If the entropy $\log s$ is positive (which happens if $f_0$, $f_1$ are expanding), then Theorem \ref{thm:ss0} says that the kneading determinant $D$ vanishes at $t=1/s$ and it does not vanish in the disk $|t|<1/s$. From the previous discussions (see matrix \ref{eq:matrix}) it is not difficult to get that
\[D=\frac{-1}{1-t} N_{2,1}.\]
Thus $t=1/s$ is the first zero of $N_{2,1}$ in the interval $(0,1)$. We conclude that $t=1/s$ is the smallest solution $t\in (0,1)$ (in fact in the open complex unit disk) to the equation
\begin{equation}\label{eq:D}
\sum_{i\geq 0} \alpha_i t^i=\sum_{i\geq 0} \beta_i t^i.
\end{equation}
Our results also imply that $\F$ in $[0,1]$ is semiconjugate to a, possibly degenerate, system $U=\{U_0,U_1\}$ in $[0,1]$ such that each $U_i$ is affine with slope $s$. Note that if $\text{Dom}(U_0)=\{0\}$ then $U_1([0,1])\subset [0,1]$ which would contradict the fact that $s>1$. Similarly one sees that $\text{Dom}(U_1)\neq \{1\}$, thus $U$ is non degenerate.

Now assume that we also have $f_0(0)=0$, $f_1(1)=1$. The semiconjugacy condition implies $U_0(0)=0$, $U_1(1)=1$ and $U_0, \, U_1$ are of the form 
\begin{align*}
U_0(x)=sx,\,\, U_1(x)=s(x-1)+1,
\end{align*}
where ${\rm Dom}(U_0)=[0,p]$, ${\rm Dom}(U_1)=[p,1]$, with $p =\varphi(q)$. To completely determine $U$ it remains to find $p$. Since $q\in \text{Dom}(f_{\alpha_n}\circ\cdots\circ f_{\alpha_0})$ for all $n$, Lemma \ref{lem:semiG} implies that 
\[
z_n := U_{\alpha_n}\circ\cdots\circ U_{\alpha_0}(p)=\varphi(f_{\alpha_n}\circ\cdots\circ f_{\alpha_0}(q)) \in [0,1].
\]
From this we get
\[
p =U_{\alpha_0}^{-1}\circ\cdots\circ U_{\alpha_n}^{-1}(z_n)=\frac{z_n}{s^{n+1}}+(1-s^{-1}) \sum_{i=0}^n \alpha_i s^{-i}\, \textrm{ for all } n.
\]
Making $n$ go to infinity gives
\[
p=(1-s^{-1}) \sum_{i=0}^{\infty} \alpha_i s^{-i}.
\]
\end{proof}

Notice that our theorem implies part (2) of Statement \ref{BHV1}
and for part (1) we only get existence of a semiconjugancy, not a
conjugacy. Note also that we did not assume $f_0$, $f_1$ expanding,
but assumed that $\F$ has positive entropy. 


\bibliography{KT_IFS}
\bibliographystyle{plain}

\end{document}